\documentclass[12pt,letter]{article}
\usepackage{amsmath,amssymb,amsthm,amstext}
\usepackage{enumerate}
\usepackage{graphicx}
\usepackage{epsfig}
\usepackage{epstopdf}
\usepackage{color}
\usepackage[margin=1 in]{geometry}
\usepackage{setspace}

\newtheorem{theorem}{Theorem}[section]
\newtheorem{remark}{Remark}[section]

\newtheorem{lemma}[theorem]{Lemma}
\newtheorem{proposition}[theorem]{Proposition}

\theoremstyle{remark}

\title{Positive solutions for fractional-order boundary value problems with or without dependence of integer-order ones}
\author{Inbo Sim\thanks{Department of Mathematics, University of Ulsan, Ulsan 44610, Republic of Korea, e-mail: ibsim@ulsan.ac.kr}~ and 
Satoshi Tanaka\thanks{Mathematical Institute, Tohoku University, Aoba 6-3, Aramaki, Aoba-ku, Sendai 980-8578, Japan, e-mail: satoshi.tanaka.d4@tohoku.ac.jp}}
\begin{document}
\small
\date{\small Dedicated to Professor Y\={u}ki Naito on the occasion of his 60th birthday}
\maketitle

\begin{abstract}
We investigate the existence, non-existence, uniqueness, and multiplicity of positive solutions to the following problem:
\begin{align}\label{P}
 \left\{
  \begin{array}{l}
   D_{0+}^\alpha u + h(t)f(u) = 0, \quad 0<t<1, \\[1ex]
   u(0)=u(1)=0,
  \end{array}
 \right.
\end{align}
where $D_{0+}^\alpha$ is the Riemann-Liouville fractional derivative of 
order $\alpha\in(1,2]$. 

Firstly, by considering the first eigenvalue $\lambda_1(\alpha)$ of the corresponding eigenvalue problem, we establish the existence of positive solutions for both sublinear and superlinear cases involving $\lambda_1(\alpha)$, thereby extending existing results in the literature. In addition, we address the issue of non-existence, which reinforces the sharpness of both hypotheses. 

Secondly, we demonstrate the uniqueness of positive solutions. For the sublinear case, we impose certain monotonicity conditions on $f$. For the superlinear case, we assume that $h$ satisfies a specific condition to ensure the uniqueness of positive solutions when $\alpha =2.$ Near $\alpha =2,$
 we prove uniqueness by leveraging the non-degeneracy of the unique solution, which represents a novel approach to studying fractional-order differential equations.

Finally, we apply this methodology to establish the multiple existence of at least three positive solutions for H\'{e}non-type problems, which is also a new contribution.

\end{abstract}


\smallskip

\setcounter{equation}{0}

\section{Introduction and Main Results} \label{sec:1}

\setcounter{section}{1} \setcounter{equation}{0} 

For fractional differential equations involving the Riemann-Liouville or Caputo fractional derivative, certain classical methods such as integration by parts and the shooting method are restrictive due to the lack of a variational structure and sufficient differentiability. Efforts have been made to develop alternatives, such as modified integration by parts formulas and variational structures suited to fractional settings \cite{KOM}. Unfortunately, these approaches are not well-suited for application to the problems currently under consideration.

Consequently, there have been relatively few results on the existence of solutions for fractional-order boundary value problems. Moreover, only a limited class of weight functions and nonlinearities can be addressed, as the available methods such as fixed-point theory \cite{JY}, \cite{LL2} and bifurcation theory \cite{Liu} are constrained in their applicability. One of the primary goals is to extend the class of weight functions and nonlinearities to include those that closely resemble those used in integer-order equations.

In this paper, we address the existence, non-existence, uniqueness, and multiplicity of positive solutions to the following problem:
\begin{align}\label{P}
 \left\{
  \begin{array}{l}
   D_{0+}^\alpha u + h(t)f(u) = 0, \quad 0<t<1, \\[1ex]
   u(0)=u(1)=0,
  \end{array}
 \right.
\end{align}
where $D_{0+}^\alpha$ is the Riemann-Liouville fractional derivative of 
order $\alpha\in(1,2]$, that is,
\begin{align*}
 D_{0+}^\alpha u(t) = \left\{
  \begin{array}{ll}
  \displaystyle\frac{1}{\Gamma(2-\alpha)} \frac{d^2}{dt^2}
    \int_0^t \frac{u(s)}{(t-s)^{\alpha-1}} ds, & 1<\alpha<2, \\[3ex]
   u''(t), & \alpha=2.
  \end{array}
 \right.
\end{align*}
When $\alpha=2$, problem \eqref{P} is 
\begin{align}\label{P2}
 \left\{
  \begin{array}{l}
   u'' + h(t)f(u) = 0, \quad 0<t<1, \\[1ex]
   u(0)=u(1)=0.
  \end{array}
 \right.
\end{align}
Now, we assume the following conditions (H) and (F):
\begin{enumerate}
 \item[(H)] $h\in C[0,1]$, $h(t)\ge 0$ for $t\in[0,1]$, 
	        and $h(t)>0$ almost everywhere on $[0,1]$;
 \item[(F)] $f\in C[0,\infty)$ and $f(s)>0$ for $s>0$.
\end{enumerate}

The first aim of this paper is to establish the existence of positive solutions for problems with a broader class of nonlinearities than those considered in existing studies, which primarily rely on fixed-point arguments. Many previous works assume that the limits $\lim_{s \to 0^+} f(s)/s$ and 
$\lim_{s \to \infty} f(s)/s$
exist and are either  $0, \infty$ or sometimes finite (see \cite{HG}, \cite{LL2}). 
In contrast, we assume the weaker conditions. 
Additionally, we consider cases where the values of $f(s)/s$ near $s=0$ and 
$s=\infty$ are separated by   
$\lambda_1(\alpha),$ the first eigenvalue of the following eigenvalue problem:
\begin{align}\label{EVP}
 \left\{
  \begin{array}{l}
   D_{0+}^{\alpha} \phi + \lambda h(t) \phi = 0, \quad 0<t<1, \\[1ex]
   \phi(0)=\phi(1)=0,
  \end{array}
 \right.
\end{align}
where $\lambda\ge 0$. 
This eigenvalue problem was briefly mentioned in \cite{HG}, \cite{Liu} in the context of applying bifurcation theory, but here, we provide a detailed proof for the benefit of the reader. A comprehensive study of eigenvalue problems with Dirichlet boundary conditions, showing a sequence of eigenvalues with their asymptotic behavior, can be found in \cite{K}. In Section \ref{sec:4}, we will prove the following result.

\begin{proposition}\label{1stEV}
 Let $\alpha\in(1,2]$.
 Assume that \textup{(H)} holds.
 Then there exists $\lambda_1(\alpha)>0$ such that problem \eqref{EVP} has a 
 positive solution when $\lambda=\lambda_1(\alpha)$ and
 problem \eqref{EVP} has no positive solution when $\lambda\ne\lambda_1(\alpha)$.
 Moreover, if $0\le\lambda<\lambda_1(\alpha)$, then $\lambda$ is not an
 eigenvalue of problem \eqref{EVP}.
\end{proposition}

As mentioned above, we will establish the existence of positive solutions for both sublinear and superlinear cases involving  $\lambda_1(\alpha),$ thereby extending existing results in the literature. The following assumptions, previously considered in \cite{BDGK}, \cite{dM}, \cite{Er}, \cite{Ko}, \cite{NT}, \cite{NT1}, were used to study the existence of multiple solutions for integer-order problems.

\begin{theorem}\label{existence}
 Let $\alpha\in(1,2]$.
 Assume that \textup{(H)} and \textup{(F)} hold and either
 \begin{equation}
  \frac{f(s)}{s} > \lambda_1(\alpha), \quad s\in(0,s_0) \quad 
   \mbox{for\ some\ } s_0>0 \quad \mbox{and} \quad
  \limsup_{s\to\infty} \frac{f(s)}{s} < \lambda_1(\alpha)
  \label{f0>lam1>foo}
 \end{equation}
 or
 \begin{equation}
  \frac{f(s)}{s} < \lambda_1(\alpha), \quad s\in(0,s_0) \quad 
   \mbox{for\ some\ } s_0>0 \quad \mbox{and} \quad
  \liminf_{s\to\infty} \frac{f(s)}{s} > \lambda_1(\alpha).
  \label{f0<lam1<foo}
 \end{equation}
 Then problem \eqref{P} has at least one positive solution.
\end{theorem}

Applying Theorem 2.14 and the proof of Corollary 3.11 
in Benmeza\"{\i}, Chentout and Henderson \cite{BCH}, 
we conclude that problem \eqref{P} has at least one positive solution 
if \textup{(H)} and \textup{(F)} hold and either
\begin{equation*}
 \liminf_{s\to0^+} \frac{f(s)}{s} > \lambda_1(\alpha) > 
 \limsup_{s\to\infty} \frac{f(s)}{s}
 \quad \mbox{or} \quad
 \limsup_{s\to0^+} \frac{f(s)}{s} < \lambda_1(\alpha) < 
 \liminf_{s\to\infty} \frac{f(s)}{s}.
\end{equation*}

Theorem \ref{existence} is sharp in the following sense.

\begin{theorem}\label{nonexistence}
 Let $\alpha\in(1,2]$.
 Assume that \textup{(H)} and \textup{(F)} hold and either
 \begin{equation*}
  \frac{f(s)}{s} > \lambda_1(\alpha) \quad \mbox{for} \ s>0~~
 \mbox{or}~~
  \frac{f(s)}{s} < \lambda_1(\alpha) \quad \mbox{for} \ s>0.
 \end{equation*}
 Then problem \eqref{P} has no positive solution.
\end{theorem}

From Proposition 2.13 and the proof of Corollary 3.10
in Benmeza\"{\i}, Chentout and Henderson \cite{BCH}, it follows that
problem \eqref{P} has no positive solution if 
\textup{(H)} and \textup{(F)} hold and either
\begin{equation*}
 \inf_{s>0} \frac{f(s)}{s} > \lambda_1(\alpha)
 \quad \mbox{or} \quad
 \sup_{s>0} \frac{f(s)}{s} < \lambda_1(\alpha).
\end{equation*}

The second aim is to study the uniqueness in ODEs, which serves as a foundation for exploring important properties such as stability and continuity. While much research has been conducted on integer-order differential equations, fewer studies have addressed the challenges posed by the lack of differentiability in solutions, as previously noted. Non-autonomous problems, in particular, are more challenging to study than autonomous ones due to technical complexities. Some results in \cite{BB}, \cite{GKKW}, \cite{ZH}, derived using fixed-point theory, require the nonlinear term $f$ to satisfy a Lipschitz condition. Our goal is to relax this restriction in the sublinear case.

\begin{theorem}\label{uniquenesssub}
 Let $\alpha\in(1,2]$.
 Assume that \textup{(H)} and \textup{(F)} hold and that 
 $f(s)$ is strictly increasing, $f(s)/s$ is strictly decreasing in $s>0$, and 
 \begin{align}\label{limf(s)/s}
  \lim_{s\to0^+} \frac{f(s)}{s} > \lambda_1(\alpha) >
  \lim_{s\to\infty} \frac{f(s)}{s}.
 \end{align}
 Then problem \eqref{P} has a unique positive solution.
\end{theorem}

A key contribution of this study is the establishment of the uniqueness of positive solutions in the superlinear case, a result that, to the best of the authors' knowledge, has not been previously achieved.

\begin{remark}\label{remoff(s)/s}
 If \textup{(F)} holds and $\displaystyle \frac{f(s)}{s}$ is strictly decreasing in $s>0$, 
 then $\displaystyle  \lim_{s\to0^+} \frac {f(s)}{s}$ exists in $(0,\infty]$ and
 $\displaystyle  \lim_{s\to\infty} \frac {f(s)}{s}$ exists in $[0,\infty)$.
\end{remark}

\begin{theorem}\label{uniquenesssuper}
 Suppose that $h\in C^1[0,1]$, $h(t)>0$ for $t\in[0,1]$ and
 \begin{align}
  -\frac{2}{t} \le \frac{h'(t)}{h(t)} \le \frac{2}{1-t}, \quad t \in (0,1).
  \label{<h'/h<}
 \end{align}
 Assume, moreover, that 
 $f\in C^1[0,\infty)$, $f(s)>0$ for $s>0$,   
 $f'(0)<\lambda_1(2)<\lim_{s\to\infty}f(s)/s$,
 $\limsup_{s\to\infty}s^{-p}f(s)<\infty$ for some $p>0$, and
 \begin{align}
  f'(s)>\frac{f(s)}{s}, \quad s>0.
  \label{f'>f/s} 
 \end{align}
 Then there exists $\alpha_1\in(1,2)$ such that problem \eqref{P} has a unique 
 positive solution for each fixed $\alpha\in(\alpha_1,2]$.
\end{theorem}

\begin{remark}\label{remoff}
 We assume that $f\in C^1[0,\infty)$ and \eqref{f'>f/s} holds. 
 Then $f(0)=0$ and $(f(s)/s)'>0$ for $s>0$.
 Hence, 
 \begin{align*}
  f'(0)=\lim_{s\to0^+} \tfrac{f(s)}{s} = \limsup_{s\to0^+} \displaystyle \frac{f(s)}{s}
 \end{align*}
 and, since $f(s)/s$ is increasing, we find that 
 $\lim_{s\to\infty}f(s)/s$ exists in $(0,\infty]$. 
 In \cite{Be}, Berestycki assumed that the ratio $f(s)/s$ is increasing to 
 obtain uniqueness results for a problem when $\alpha=2$ and $h(t)$ is constant.
\end{remark}

The Lipschitz condition ensures the uniqueness of solutions
(not restricted to positive solutions).
For example, applying Theorem 2 in Benmezai and N. Benkaci-Ali \cite{BB}, 
we find that problem \eqref{P} has a unique solution 
if there exists $c\in[0,\lambda_1(\alpha))$ such that
\begin{equation*}
 |f(x)-f(y)| \le c |x-y|, \quad x, y \in [0,\infty).
\end{equation*}
When (F) holds and $f(0)=0$, $u=0$ is the unique solution of \eqref{P} and
problem \eqref{P} has no positive solution.
In this case, we have $0<f(s)<cs<\lambda_1(\alpha)s$ for $s>0$, and then
the latter condition in Theorem 1.3 holds.

The final aim is to establish a multiple existence result, specifically demonstrating the existence of at least three positive solutions for fractional-order differential equations. While some research has been conducted in this area, it has primarily focused on the shape of the nonlinearity $f$, as the methods are often limited to fixed-point theories (see \cite{WW} and references therein). In contrast, we present a concrete example that focuses on the weight function for the superlinear nonlinear term to demonstrate the existence of at least three positive solutions. As far as the authors are aware, this is the first result of its kind achieved using approaches similar to those employed in the proof of uniqueness for the superlinear case.

\begin{theorem}\label{nonuniquenessexample}
 Let $l>1$ and $p>1$ satisfy $(p-1)l\ge 4$.
 Then there exist $\delta_0>0$ such that, for almost everywhere 
 $\delta\in(-\delta_0,\delta_0)$, there exists $\alpha_1\in(1,2)$ 
 for which  
 \begin{align}
  \left\{
  \begin{array}{l}
   D_{0+}^\alpha u + \left| t -\dfrac{1}{2} + \delta \right|^l u^p = 0, 
   \quad 0<t<1, \\[3ex]
   u(0)=u(1)=0
  \end{array} 
  \right.
  \label{fracHenon}
 \end{align}
 has at least three positive solutions for $\alpha\in(\alpha_1,2]$.
\end{theorem} 

This paper is organized as follows: In Section 2, we discuss the properties of the Green's function associated with the fractional operator under Dirichlet boundary conditions and present several solution estimates. Sections 3 and 4 explore the order topology, define an operator on a suitable space equipped with this topology, and establish the existence, estimates, and continuity of the first eigenvalue for \eqref{EVP}. In Section 5, we prove the non-existence of positive solutions for \eqref{1stEV}.
Sections 6 and 7 focus on demonstrating the existence of positive solutions for the sublinear and superlinear cases using Schauder fixed point theory and Leray-Schauder degree theory, respectively. In Sections 8 and 9, we establish the uniqueness of positive solutions for the sublinear and superlinear cases. Finally, the last two sections are devoted to proving the existence of at least three positive solutions for a special form of \eqref{P}, specifically the H\'{e}non-type equation.

We will use the following notations:
\begin{align*}
 & ||u||_\infty := \max_{t\in[0,1]} |u(t)|; \qquad 
   ||u||_{C_{2-\alpha}} := \max_{t\in[0,1]}t^{2-\alpha}|u(t)|; \\
 & C_{2-\alpha}^1[0,1] := \{ u \in C[0,1] \cap C^1(0,1] : 
   t^{2-\alpha}u'(t) \in C[0,1], u(0)=u(1)=0 \}.
\end{align*}

\section{Green's function} \label{sec:2}

\setcounter{section}{2} \setcounter{equation}{0} 

In this section, we discuss the properties of the Green's function associated with the fractional operator under Dirichlet boundary conditions and present several solution estimates that will be frequently used.  Bai and L\"{u} \cite[Lemma 2.3]{BL} introduced the Green's function of 
$-D_{0+}^\alpha$ with the Dirichlet data as follows.

\begin{proposition}\label{Green}
 Let $\alpha\in(1,2]$.
 Then for a given $g\in C[0,1]$, the unique solution of
 \begin{align}\label{Gg}
 \left\{
  \begin{array}{l}
   -D_{0+}^\alpha v = g(t), \quad 0<t<1, \\[1ex]
   v(0)=v(1)=0
  \end{array}
 \right.
 \end{align}
 is represented by
 \begin{align*}
  v(t) = \int_0^1 G(t,s,\alpha) g(s) ds, \quad t\in[0,1],
 \end{align*}
 where
 \begin{align*}
  G(t,s,\alpha) = \left\{
   \begin{array}{ll}
    \dfrac{(t(1-s))^{\alpha-1}-(t-s)^{\alpha-1}}{\Gamma(\alpha)}, &
     0 \le s \le t \le 1, \\[2ex]
    \dfrac{(t(1-s))^{\alpha-1}}{\Gamma(\alpha)}, &
     0 \le t \le s \le 1.
   \end{array}
  \right.
 \end{align*}
\end{proposition}

We have the following properties of Green's function $G(t,s,\alpha)$,
which obtained in \cite[Lemma 2.4 and its proof]{BL} and 
\cite[Theorems 1.1 and 1.2]{JY}.

\begin{lemma}\label{propofG}
 Let $\alpha\in(1,2]$.
 Then Green's function $G(t,s,\alpha)$ satisfies the following 
 \textup{(i)--(iv):}
 \begin{enumerate}
  \item[\textup{(i)}] $G\in C([0,1]\times[0,1]\times(1,2])$, $G(t,s,\alpha)>0$ 
	for $t$, $s \in(0,1)$, 
	$G(0,s,\alpha)=G(1,s,\alpha)=0$ for $s\in[0,1]$, and
	$G(t,0,\alpha)=G(t,1,\alpha)=0$ for $t\in[0,1]$\textup{;}
  \item[\textup{(ii)}] $G(t,s,\alpha)$ is decreasing with respect to $t$ for $s \le t$ and 
	increasing with respect to $t$ for $t \le s$, and hence
        \begin{align*}
	  \max_{t\in[0,1]} G(t,s,\alpha)
           = G(s,s,\alpha)
           = \frac{1}{\Gamma(\alpha)} (s(1-s))^{\alpha-1}
           \quad \mbox{for} \ s \in [0,1]\textup{;}
	\end{align*}
  \item[\textup{(iii)}] for $t$, $s \in [0,1]$,
	\begin{align*}
	 \frac{\alpha-1}{\Gamma(\alpha)} t^{\alpha-1}(1-t) s(1-s)^{\alpha-1}
	 \le G(t,s,\alpha) \le 
         \frac{1}{\Gamma(\alpha)} t^{\alpha-1}(1-t) (1-s)^{\alpha-2}\textup{;}
	\end{align*}
  \item[\textup{(iv)}] for $t$, $s \in [0,1]$,
	\begin{align*}
	 \frac{\alpha-1}{\Gamma(\alpha)} t(1-t) s(1-s)^{\alpha-1}
	 \le t^{2-\alpha} G(t,s,\alpha) \le 
         \frac{1}{\Gamma(\alpha)} s(1-s)^{\alpha-1}.
	\end{align*}
 \end{enumerate}
\end{lemma}

\begin{lemma}\label{lemv>}
 Let $\alpha\in(1,2]$ and let $g\in C[0,1]$ satisfy $g(t)\ge 0$ for $t\in[0,1]$.
 If $v$ be a solution of \eqref{Gg}, then $v$ satisfies
 \begin{align}\label{v>}
  t^{2-\alpha} v(t) 
  \ge (\alpha-1) t(1-t) ||v||_{C_{2-\alpha}}, \quad t \in [0,1].
 \end{align}
\end{lemma}

\begin{proof} 
 This result has been proved in the proof of Lemma 4.1 in \cite{JY}.
 However, we give a proof for readers' convenience.
 From (iv) of Lemma \ref{propofG}, it follows that 
 \begin{align*}
  t^{2-\alpha} v(t) = t^{2-\alpha} \int_0^1 G(t,s,\alpha) g(s) ds
  \ge \frac{\alpha-1}{\Gamma(\alpha)} t(1-t) 
        \int_0^1 s(1-s)^{\alpha-1} g(s) ds
 \end{align*}
 and 
 \begin{align*}
  t^{2-\alpha} v(t) 
   \le \frac{1}{\Gamma(\alpha)} \int_0^1 s(1-s)^{\alpha-1} g(s) ds,
 \end{align*}
 which imply 
 \begin{align*}
  ||v||_{C_{2-\alpha}} 
  \le \frac{1}{\Gamma(\alpha)} \int_0^1 s(1-s)^{\alpha-1} g(s) ds.
 \end{align*}
 Combining these inequalities, we have \eqref{v>}.
\end{proof} 

\begin{lemma}\label{|v'-v'|}
 Let $\alpha\in(1,2]$, $g\in C[0,1]$ and $v$ be a unique solution of 
 \eqref{Gg}. 
 Then
 \begin{align*}
 \left| t_1^{2-\alpha}v'(t_1) - t_2^{2-\alpha}v'(t_2) \right|
   \le \frac{||g||_\infty}{\Gamma(\alpha)} 
        |t_2^{2-\alpha} - t_1^{2-\alpha}|
     + \frac{2 ||g||_\infty}{\Gamma(\alpha)}|t_2-t_1|^{\alpha-1}
 \end{align*}
 for $0\le t_1 \le t_2 \le 1$.
\end{lemma}

\begin{proof} 
 Since
 \begin{align*}
  v(t) = \int_0^1 \dfrac{(t(1-s))^{\alpha-1}}{\Gamma(\alpha)} g(s) ds
       - \int_0^t \dfrac{(t-s)^{\alpha-1}}{\Gamma(\alpha)} g(s) ds,
 \end{align*}
 it follows that
 \begin{align}\label{v'}
 v'(t) = \frac{\alpha-1}{\Gamma(\alpha)} t^{\alpha-2} 
      \int_0^1 (1-s)^{\alpha-1} g(s) ds
    - \frac{\alpha-1}{\Gamma(\alpha)} \int_0^t (t-s)^{\alpha-2} g(s) ds
 \end{align}
for $t\in(0,1]$.
Assume that $0\le t_1 \le t_2 \le 1$. 
Then we have
\begin{align*}
 \frac{\Gamma(\alpha)}{\alpha-1} 
 \Bigl( t_1^{2-\alpha}&v'(t_1) - t_2^{2-\alpha}v'(t_2) \Bigr) \\
  & = - t_1^{2-\alpha} \int_0^{t_1} (t_1-s)^{\alpha-2} g(s) ds 
      + t_2^{2-\alpha} \int_0^{t_2} (t_2-s)^{\alpha-2} g(s) ds \\
  & = (t_2^{2-\alpha} - t_1^{2-\alpha}) 
       \int_0^{t_1} (t_1-s)^{\alpha-2} g(s) ds \\
  & \quad \ + t_2^{2-\alpha} \left( 
          \int_0^{t_2} (t_2-s)^{\alpha-2} g(s) ds
        - \int_0^{t_1} (t_1-s)^{\alpha-2} g(s) ds 
          \right),
\end{align*}
and hence
\begin{multline*}
  \frac{\Gamma(\alpha)}{\alpha-1} 
 \left| t_1^{2-\alpha}v'(t_1) - t_2^{2-\alpha}v'(t_2) \right| \\
   \le |t_2^{2-\alpha} - t_1^{2-\alpha}| 
        \frac{t_1^{\alpha-1}}{\alpha-1} ||g||_\infty
     + \left| \int_0^{t_2} (t_2-s)^{\alpha-2} g(s) ds
        - \int_0^{t_1} (t_1-s)^{\alpha-2} g(s) ds \right|.
\end{multline*}
We observe that
\begin{align*}
 & \left| \int_0^{t_2} (t_2-s)^{\alpha-2} g(s) ds
 - \int_0^{t_1} (t_1-s)^{\alpha-2} g(s) ds \right| \\
 & = \left| \int_{t_1}^{t_2} (t_2-s)^{\alpha-2} g(s) ds 
 + \int_0^{t_1} [(t_2-s)^{\alpha-2}-(t_1-s)^{\alpha-2}] g(s) ds \right| \\
 & \le ||g||_\infty \int_{t_1}^{t_2} (t_2-s)^{\alpha-2} ds
   + ||g||_\infty \int_0^{t_1} [(t_1-s)^{\alpha-2}-(t_2-s)^{\alpha-2}] ds \\ 
 & = \frac{||g||_\infty}{\alpha-1}(t_2-t_1)^{\alpha-1} 
     +\frac{||g||_\infty}{\alpha-1} \left((t_2-t_1)^{\alpha-1}
                               - (t_2^{\alpha-1}-t_1^{\alpha-1})\right) \\
 & \le \frac{2 ||g||_\infty}{\alpha-1}(t_2-t_1)^{\alpha-1}.
\end{align*}
Therefore we have
\begin{align*}
 \frac{\Gamma(\alpha)}{\alpha-1} 
 \left| t_1^{2-\alpha}v'(t_1) - t_2^{2-\alpha}v'(t_2) \right|
   \le \frac{||g||_\infty}{\alpha-1} |t_2^{2-\alpha} - t_1^{2-\alpha}|
     + \frac{2 ||g||_\infty}{\alpha-1}|t_2-t_1|^{\alpha-1}.
\end{align*}
\end{proof} 

\section{Order topology} \label{sec:3}

\setcounter{section}{3} \setcounter{equation}{0} 

In this section, we examine the order topology to address our problems and define an operator on a suitable space equipped with this topology. First, we prove the following result.

\begin{lemma}\label{|y|<}
 Let $\alpha\in(1,2]$ and $y\in C_{2-\alpha}^1[0,1]$.
 Then
 \begin{align}
  |y(t)| \le \frac{||y'||_{C_{2-\alpha}}}{\alpha-1}
   \min\{t^{\alpha-1},1-t^{\alpha-1}\}, \quad t \in [0,1].
  \label{|y|<1}
 \end{align}
 In particular,
 \begin{align}
  |y(t)| \le \frac{2||y'||_{C_{2-\alpha}}}{\alpha-1} t^{\alpha-1}(1-t), 
   \quad t \in [0,1].
  \label{|y|<2}
 \end{align}
\end{lemma}

\begin{proof} 
 Inequalities \eqref{|y|<1} and \eqref{|y|<2} hold at $t=0$, $1$.
 Let $t\in(0,1)$.
 We observe that
 \begin{equation*}
  |y(t)|   
    = \left| \int_t^1 s^{\alpha-2} s^{2-\alpha} y'(s) ds \right|
   \le ||y'||_{C_{2-\alpha}} \int_t^1 s^{\alpha-2} ds \\
  = \frac{||y'||_{C_{2-\alpha}}}{\alpha-1} (1-t^{\alpha-1}). 
 \end{equation*}
 For each $\varepsilon\in(0,t)$, we have
 \begin{align*}
  |y(t)|  = \left| y(\varepsilon) + \int_\varepsilon^t y'(s) ds \right|
  & \le |y(\varepsilon)| 
      + ||y'||_{C_{2-\alpha}} \int_\varepsilon^t s^{\alpha-2} ds \\
  & \le |y(\varepsilon)|
      + \frac{||y'||_{C_{2-\alpha}}}{\alpha-1}
         (t^{\alpha-1}-\varepsilon^{\alpha-1}). 
 \end{align*}
 Letting $\varepsilon\to0^+$, we get
 \begin{align*}
  |y(t)| \le \frac{||y'||_{C_{2-\alpha}}}{\alpha-1} t^{\alpha-1}.
 \end{align*}
 Therefore, \eqref{|y|<1} holds.

 If $0<t\le 2^{-1/(\alpha-1)}$, then
 \begin{align*}
  2 t^{\alpha-1}(1-t) \ge 2 t^{\alpha-1}(1-2^{-1/(\alpha-1)})
  \ge t^{\alpha-1} 
  \ge \min \{t^{\alpha-1},1-t^{\alpha-1}\}.
 \end{align*}
 If $2^{-1/(\alpha-1)}\le t <1$, then
 \begin{align*}
  2 t^{\alpha-1}(1-t) \ge 2 (2^{-1/(\alpha-1)})^{\alpha-1}(1-t)
  \ge 1-t^{\alpha-1} \ge \min \{t^{\alpha-1},1-t^{\alpha-1}\}.
 \end{align*}
 Consequently, \eqref{|y|<2} holds.
\end{proof} 

Next, let $\alpha\in(1,2]$.
According to Deimling \cite[Section 19.6]{De}, we consider the order topology
as follows.
See also Amann \cite{Am}.
Let $X=C[0,1]$ be the Banach space equipped with the uniform norm
$||\,\cdot\,||_\infty$.
We set
\begin{gather*}
 K = \{ x \in X : x(t) \ge 0 \ \mbox{for} \ t \in [0,1] \}; \qquad
 e(t) = t^{\alpha-1}(1-t) \quad \textup{for} \ t \in [0,1].
\end{gather*}
Then $K$ is a cone.
We define
\begin{gather*}
 X_e= \bigcup_{\rho>0} \{ x \in X :
  -\rho e(t) \le x(t) \le \rho e(t) \ \mbox{for} \ t \in [0,1] \}, \\
 ||x||_e  = \inf \{\rho>0 :
   -\rho e(t) \le x(t) \le \rho e(t) \ \mbox{for} \ t \in [0,1] \} \
 \mbox{for} \ w \in X_e.
\end{gather*}
Then $X_e$ is a subspace of $X$ and $||\cdot||_e$ is a norm on $X_e$.
If $x \in X_e$, then $|x(t)|\le ||x||_e e(t)$ for $t\in[0,1]$, and hence
\begin{align}\label{||x||<||x||e}
 ||x||_\infty \le ||x||_e, \quad x \in X_e.
\end{align}
Moreover, Lemma \ref{|y|<} implies $C_{2-\alpha}^1[0,1]\subset X$ and
\begin{align}\label{||x||e<||x||C2-a}
 ||x||_e \le \frac{2}{\alpha-1} ||x'||_{C_{2-\alpha}}, 
 \quad x \in C_{2-\alpha}^1[0,1].
\end{align}
Using \eqref{||x||<||x||e}, we can show that $(X_e,||\cdot||_e)$ is complete, hence a Banach space. 

We set
\begin{align*}
  K_e = X_e \cap K = \{ x \in K :
                       0 \le x(t) \le \rho e(t) \ \mbox{for\ some}\ \rho>0 \}.
\end{align*}
We conclude that $K_e$ is a cone in $X_e$ and 
\begin{align*}
 \mbox{int}(K_e) = \{ x \in K_e : 
                       x(t) \ge \delta e(t) \ \mbox{for\ some}\ \delta>0 \},
\end{align*}
where $\textup{int}(K_e)$ is the interior of $K_e$.

\begin{lemma}\label{compact}
 Let $\alpha\in(1,2]$.
 Assume that \textup{(H)} holds.
 Then the operator
 \begin{align*}
  (Tx)(t):= \int_0^1 G(t,s,\alpha) h(s) x(s) ds
 \end{align*}
 satisfies $T\in\mathcal{L}(X_e)$, 
 $T(K_e\setminus\{0\})\subset\textup{int}\,(K_e)$ and 
 $T$: $X_e\to X_e$ is compact.
\end{lemma}

\begin{proof} 
 Since
 \begin{align}
  \int_0^1 G(t,s,\alpha) ds = \frac{1}{\alpha\Gamma(\alpha)} e(t) 
  = \frac{1}{\Gamma(\alpha+1)} e(t), \quad t \in [0,1],
  \label{intG}
 \end{align}
 we get
 \begin{align}\label{<Tx<}
  - \frac{1}{\Gamma(\alpha+1)} ||h||_\infty ||x||_\infty e(t)
  \le (Tx)(t) \le  
    \frac{1}{\Gamma(\alpha+1)} ||h||_\infty ||x||_\infty e(t)
 \end{align}
 for $x\in X_e$, $t \in [0,1]$.
 Hence, $T$ maps $X_e$ into itself.
 Since
 \begin{align*}
  ||Tx||_e \le \frac{1}{\Gamma(\alpha+1)} ||h||_\infty ||x||_\infty
  \le \frac{1}{\Gamma(\alpha+1)} ||h||_\infty ||x||_e,
  \quad x \in X_e,
 \end{align*}
 by \eqref{||x||<||x||e}, we see that $T\in\mathcal{L}(X_e)$.

 Let $x \in K_e\setminus\{0\}$.
 Since $x(t)\ge0$ and $x(t)\not\equiv0$ on $[0,1]$, we have 
 $(Tx)(t)\ge0$ and $(Tx)(t)\not\equiv0$ on $[0,1]$ and hence 
 $||Tx||_{C_{2-\alpha}}>0$, because of Lemma \ref{propofG}  (iii).
 Lemma \ref{lemv>} implies that 
 $(Tx)(t)\ge (\alpha-1)||Tx||_{C_{2-\alpha}}e(t)$ for $t\in[0,1]$.
 Hence, we find that $Tx\in \mbox{int}(K_e)$. 

 Finally we prove that $T$ : $X_e\to X_e$ is compact.
 Let $\{x_n\}_{n=1}^\infty$ be a sequence in $X_e$ such that
 $||x_n||_e \le M$ for some $M>0$.
 By \eqref{||x||<||x||e}, we have $||x_n||_\infty \le M$.
 We note that $T$ : $X\to X$ is compact.
 Thus, $\{Tx_n\}_{n=1}^\infty$ contains a subsequences converging in $X$.
 Lemma \ref{|v'-v'|} ensures that $\{t^{2-\alpha}(Tx_n)'(t)\}_{n=1}^\infty$ 
 is equicontinuos on $[0,1]$. 
 Using \eqref{v'}, we get
 \begin{align*}
  \frac{\Gamma(\alpha)}{\alpha-1} |t^{2-\alpha}(Tx_n)'(t)|
  & \le \int_0^1 (1-s)^{\alpha-1} h(s)|x_n(s)| ds
        + t^{2-\alpha} \int_0^t (t-s)^{\alpha-2} h(s)|x_n(s)| ds \\
  & \le ||h||_\infty M \int_0^1 (1-s)^{\alpha-1} ds
      + ||h||_\infty M t^{2-\alpha} \int_0^t (t-s)^{\alpha-2} ds \\
  & = \frac{||h||_\infty M}{\alpha} + \frac{||h||_\infty M}{\alpha-1} t,
       \quad t\in[0,1].
 \end{align*}
 Consequently, by the Arzel\`{a}-Ascoli theorem, there exist a subsequence 
 $\{x_{n_k}\}_{k=1}^\infty$ and $y$, $z\in X$ such that 
 $||y_k-y||_\infty\to0$ and $||z_k-z||_\infty\to0$ as $k\to\infty$,
 where $y_k:=Tx_{n_k}$ and $z_k(t):=t^{2-\alpha}y_k'(t)$.
 
 We will prove that $y\in C^1(0,1]$ and $z(t)=t^{2-\alpha}y'(t)$ 
 for $t \in[0,1]$.
 Since 
 \begin{align*}
  y_k(t) = -\int_t^1 s^{\alpha-2} s^{2-\alpha} y_k'(s) ds
         = -\int_t^1 s^{\alpha-2} z_k(s) ds, \quad t \in (0,1].
 \end{align*}
 Letting $k\to\infty$, we have
 \begin{align*}
  y(t) = -\int_t^1 s^{\alpha-2} z(s) ds, \quad t \in (0,1],
 \end{align*}
 which implies that $y\in C^1(0,1]$ and $t^{2-\alpha}y'(t)=z(t)$ for 
 $t\in(0,1]$.
 Since $z\in C[0,1]$, we see that
 \begin{align*}
  \lim_{t\to0^+} t^{2-\alpha}y'(t) = z(0),
 \end{align*}
 so that $t^{2-\alpha}y'(t)=z(t)$ for $t\in[0,1]$.
 Since $(Tx_n)(0)=(Tx_n)(1)=0$, we see that $y\in C_{2-\alpha}^1[0,1]$.

 By \eqref{||x||e<||x||C2-a}, we have
 \begin{align*}
  ||y_k-y||_e \le \frac{2}{\alpha-1} ||y_k'-y'||_{C_{2-\alpha}}
  = \frac{2}{\alpha-1} ||z_k-z||_\infty.
 \end{align*}
 Since $||z_k-z||_\infty\to0$, we conclude that $||y_k-y||_e\to0$ as 
 $k\to\infty$.
 Consequently, $T$: $X_e\to X_e$ is compact.
\end{proof} 

We will use the following result, 
which is a kind of the Krein-Rutman theorem (\cite{KR}).
See, for the proof, Deimling \cite[Theorem 19.3]{De}.

\begin{proposition}\label{KRThm}
 Let $X$ be a Banach space, $K\subset X$ a cone with 
 $\textup{int}\,(K)\ne\emptyset$, $T\in\mathcal{L}(X)$, compact and 
 $T(K\setminus\{0\})\subset\textup{int}\,(K)$.
 Then the following \textup{(i)--(iii)} hold\textup{:}
 \begin{enumerate}
  \item[\textup{(i)}] $r(T)>0$ and $r(T)$ is a simple eigenvalue 
	with an eigenvector $v\in\textup{int}\,(K)$ and there is no other 
	eigenvalue with an eigenvector belonging to $K\setminus\{0\}$\textup{;}
  \item[\textup{(ii)}] $|\mu|<r(T)$ for all eigenvalues $\mu\ne r(T)$\textup{;}
  \item[\textup{(iii)}] If $S\in\mathcal{L}(X)$ and $(S-T)(x)\in\textup{int}\,(K)$ for 
	$x\in K\setminus\{0\}$, then $r(S)>r(T)$,
 \end{enumerate}
 where, $r(T)$ is the spectral radius of $T$.
\end{proposition}

\section{First eigenvalue} \label{sec:4}

\setcounter{section}{4} \setcounter{equation}{0} 

In this section, we prove Proposition \ref{1stEV}, provide estimates, and demonstrate the continuity of 
$\lambda_1(\alpha)$.

\bigskip

\begin{proof}[Proof of Proposition \ref{1stEV}]
 By Proposition \ref{Green}, problem \eqref{EVP} is equivalent to 
 the equation $Tx=\lambda x$, where
 \begin{align*}
  (Tx)(t) := \int_0^1 G(t,s,\alpha) h(s) x(s) ds, \quad t\in[0,1].
 \end{align*}
 We note that every solution $x\in C[0,1]$ of $Tx=\lambda x$ satisfies
 $x\in X_e$, because of \eqref{<Tx<}.
 If $\lambda=0$, then $x=0$ is the unique solution of $Tx=\lambda x$,
 because of Proposition \ref{Green}.
 Therefore, $\lambda=0$ is not an eigenvalue of problem \eqref{EVP}.

 We use Proposition \ref{KRThm} with $X=X_e$ and $K=K_e$.
 Then we conclude that $r(T)>0$ and that 
 problem \eqref{EVP} with $\lambda=1/r(T)$ has a positive solution in 
 $\mbox{int}(K_e)$ 
 and problem \eqref{EVP} has no positive solution when $\lambda\ne1/r(T)$.
 By setting $\lambda_1(\alpha)=1/r(T)$, the proof is complete.
\end{proof} 

\begin{lemma}\label{1steigenvalue}
 Let $\alpha\in(1,2]$.
 Assume that \textup{(H)} holds.
 Then the first eigenvalue satisfies 
 \begin{equation*}
  \frac{1}{\left\| \int_0^1 G(\,\cdot\,,s,\alpha) h(s) ds \right\|_\infty} 
  \le \lambda_1(\alpha)
  \le 
  \frac{1}{(\alpha-1)\left\| \int_0^1 s^{\alpha-1}(1-s) G(\,\cdot\,,s,\alpha) h(s) ds \right\|_{C_{2-\alpha}}}.
 \end{equation*}
\end{lemma}

\begin{proof} 
 Let $\phi$ be a positive solution of \eqref{EVP}
 with $\lambda=\lambda_1(\alpha)$.
 Since
 \begin{align*}
  0 < \phi(t) = \lambda_1(\alpha) \int_0^1 G(t,s,\alpha) h(s) \phi(s) ds,
  \quad t\in (0,1),
 \end{align*}
 we get
 \begin{align*}
  0 < ||\phi||_\infty 
  \le \lambda_1(\alpha) \left\| \int_0^1 G(\,\cdot\,,s,\alpha) h(s) ds \right\|_\infty
       ||\phi||_\infty, 
 \end{align*}
 which implies
  $\lambda_1(\alpha)
   \ge 1/\| \int_0^1 G(\,\cdot\,,s,\alpha) h(s) ds \|_\infty$.
 Lemma \ref{lemv>} ensures that
 \begin{align*}
  t^{2-\alpha} \phi(t) 
  \ge (\alpha-1) t(1-t) ||\phi||_{C_{2-\alpha}}, \quad t \in [0,1].
 \end{align*}
 We observe that
 \begin{align*}
  t^{2-\alpha} \phi(t) 
   & = t^{2-\alpha} \lambda_1(\alpha) \int_0^1 G(t,s,\alpha) h(s) 
       s^{\alpha-2} s^{2-\alpha} \phi(s) ds \\
   & \ge \lambda_1(\alpha) (\alpha-1) t^{2-\alpha} 
         \int_0^1 G(t,s,\alpha) h(s) s^{\alpha-2}
          s(1-s) ds ||\phi||_{C_{2-\alpha}}
 \end{align*}
 for $t\in[0,1]$, and hence
 \begin{align*}
  ||\phi||_{C_{2-\alpha}} 
   \ge \lambda_1(\alpha) (\alpha-1)  
   \left\| \int_0^1 s^{\alpha-1}(1-s) G(\,\cdot\,,s,\alpha) h(s) ds \right\|_{C_{2-\alpha}}
   ||\phi||_{C_{2-\alpha}}.
 \end{align*}
 Then 
  $\lambda_1(\alpha) \le 1/((\alpha-1)\| \int_0^1 s^{\alpha-1}(1-s) 
   G(\,\cdot\,,s,\alpha) h(s) ds \|_{C_{2-\alpha}})$.
\end{proof} 

\begin{proposition}\label{estimate1stEV}
 Let $\alpha\in(1,2]$.
 Assume that \textup{(H)} holds.
 Then the first eigenvalue $\lambda_1(\alpha)$ satisfies
 \begin{equation*}\label{estimate1stEV}
  \frac{1}{||h||_\infty} 
  \le \frac{\alpha^\alpha \Gamma(\alpha+1)}{(\alpha-1)^{\alpha-1}||h||_\infty}
  \le \lambda_1(\alpha)
  \le \frac{4\Gamma(\alpha)}{(\alpha-1)^2\int_0^1 s^{\alpha} (1-s)^\alpha h(s) ds}.
 \end{equation*}
\end{proposition}

\begin{proof} 
 Recalling \eqref{intG}, we see that
 \begin{align*}
  \int_0^1 G(t,s,\alpha) h(s) ds 
    & \le ||h||_\infty \int_0^1 G(t,s,\alpha) ds
      = \frac{||h||_\infty}{\Gamma(\alpha+1)} t^{\alpha-1}(1-t) \\
    & \le \frac{(\alpha-1)^{\alpha-1}}{\alpha^\alpha \Gamma(\alpha+1)}
          ||h||_\infty 
      \le ||h||_\infty, \quad t \in [0,1].
 \end{align*}
 From (iii) of Lemma \ref{propofG}, it follows that
 \begin{align*}
  t^{2-\alpha} \int_0^1 s^{\alpha-1}(1-s) G(t,s,\alpha) h(s) ds
  & \ge \frac{\alpha-1}{\Gamma(\alpha)} t(1-t)
        \int_0^1 s^\alpha (1-s)^\alpha h(s) ds 
 \end{align*}
 for $t\in[0,1]$.
 Hence we have
 \begin{equation*}
  (\alpha-1)\left\| \int_0^1 s^{\alpha-1}(1-s) G(\,\cdot\,,s,\alpha) h(s) ds \right\|_{C_{2-\alpha}}
  \ge \frac{(\alpha-1)^2}{4\Gamma(\alpha)} 
        \int_0^1 s^\alpha (1-s)^\alpha h(s) ds.
 \end{equation*}
 Consequently, Lemma \ref{1steigenvalue} implies
 Proposition \ref{estimate1stEV}.
\end{proof} 

\begin{proposition}\label{contioflam1}
 The first eigenvalue $\lambda_1(\alpha)$ is continuous with respect to 
 $\alpha\in(1,2]$.
\end{proposition}

\begin{proof} 
 Assume to the contrary that there exist $\alpha_0\in(1,2]$ and   
 $\{\alpha_n\}\subset(1,2]$ such that $\alpha_n\to\alpha_0$ and 
 $\{\lambda_1(\alpha_n)\}_{n=1}^\infty$ does not converges to $\lambda_1(\alpha_0)$.
 Proposition \ref{estimate1stEV} implies that 
 $\{\lambda_1(\alpha_n)\}_{n=1}^\infty$ is bounded, 
 and hence it contains a convergent subsequence.
 We denote the convergent subsequence by $\{\lambda_1(\alpha_n)\}_{n=1}^\infty$ again and 
 its limit by $\lambda$.
 Then $\lambda\ne\lambda_1(\alpha_0)$.
 We take $\alpha_*\in(1,2)$ and $\lambda_*>0$ 
 for which $\alpha_n\in[\alpha_*,2]$ and 
 $0<\lambda_1(\alpha_n)\le\lambda_*$ for $n\ge1$.
 For each $n\ge1$, let $\phi_n$ be an eigenfunction of \eqref{EVP} 
 corresponding to $\lambda=\lambda_1(\alpha_n)$ satisfying $\phi_n(t)>0$ 
 for $t\in(0,1)$ and $||\phi_n||_\infty=1$.
 From (i) of Lemma \ref{propofG}, it follows that $G(t,s,\alpha)$ is uniformly 
 continuous on $[0,1]\times[0,1]\times[\alpha_*,2]$.
 Hence, for each $\varepsilon>0$, there exists $\delta>0$ such that
 if $t_1$, $t_2\in[0,1]$, $|t_1-t_2|<\delta$, $s\in[0,1]$, 
 $\alpha\in[\alpha_*,2]$, then
 \begin{align*}
  |G(t_1,s,\alpha)-G(t_2,s,\alpha)|< \frac{\varepsilon}{\lambda_* ||h||_\infty}.
 \end{align*}
 Therefore, if $|t_1-t_2|<\delta$, then
 \begin{align*}
  |\phi_n(t_1)-\phi_n(t_2)| 
   & = \left|\lambda_1(\alpha_n)
   \int_0^1 (G(t_1,s,\alpha_n)-G(t_2,s,\alpha_n)) h(s) \phi_n(s) ds \right| \\
   & \le \lambda_* ||h||_\infty 
     \int_0^1 |G(t_1,s,\alpha_n)-G(t_2,s,\alpha_n)| ds 
     < \varepsilon.
 \end{align*}
 This means that $\{\phi_n\}_{n=1}^\infty$ is equicontinuos on $[0,1]$.
 By the Arzel\`{a}-Ascoli theorem, there exists a subsequence $\{\phi_{n_k}\}_{k=1}^\infty$ 
 of $\{ \phi_n \}_{n=1}^\infty$ that converges to $\phi_0$ uniformly on $[0,1]$ 
 for some $\phi_0 \in C[0,1]$.
 Since
 \begin{align*}
  \phi_{n_k}(t) = \lambda_1(\alpha_{n_k}) 
   \int_0^1 G(t,s,\alpha_{n_k}) h(s) \phi_{n_k}(s) ds, \quad t\in[0,1],
 \end{align*}
 Lebesgue's dominated convergence theorem, we have
 \begin{align}
  \phi_0(t) = \lambda \int_0^1 G(t,s,\alpha_0) h(s) \phi_0(s) ds, 
  \quad t \in [0,1],
  \label{phi0=int}
 \end{align}
 which means that $\phi_0$ is a nonnegative nontrivial solution of \eqref{EVP}.
 Since $||\phi_0||_\infty=1$, from \eqref{phi0=int}, it follows that
 $\phi_0(t)>0$ for $t\in(0,1)$, and hence $\phi_0$ is a positive 
 solution of \eqref{EVP}.
 Proposition \ref{1stEV} implies $\lambda$ must be $\lambda_1(\alpha)$.
 This is a contradiction.
\end{proof} 

\section{Non-existence of positive solutions} \label{sec:5}

\setcounter{section}{5} \setcounter{equation}{0} 

In this section, we show the non-existence of positive solutions for \eqref{P} under the condition $\frac{f(s)}{s} > \lambda_1(\alpha)$ or $\frac{f(s)}{s} < \lambda_1(\alpha)$ for $s>0$.
Let $\phi_1$ denote a positive solution of \eqref{EVP} with
$\lambda=\lambda_1(\alpha)$ satisfying $||\phi_1||_\infty=1$. 

\begin{lemma}\label{v>lamintGhvhasnosol}
 Let $\alpha\in(1,2]$.
 Assume that \textup{(H)} holds and that $\overline{h}\in C[0,1]$, 
 $\overline{h}(t)\ge \lambda_1(\alpha)h(t)$ and 
 $\overline{h}(t)\not\equiv \lambda_1(\alpha)h(t)$ for $t\in[0,1]$.
 Then there is no bounded measurable function $v(t)$ on $[0,1]$ satisfying
 \begin{align}\label{v>lamintGhv}
  v(t) \ge \int_0^1 G(t,s,\alpha) \overline{h}(s) v(s) ds, \quad t \in [0,1]
 \end{align}
 and
 \begin{align}\label{v>ct(1-t)}
  v(t)\ge ct^{\alpha-1}(1-t), \quad t \in [0,1]
 \end{align}
 for some constant $c>0$.
\end{lemma}

\begin{proof} 
 We suppose that there exists a bounded measurable function $v(t)$ satisfying 
 \eqref{v>lamintGhv} and \eqref{v>ct(1-t)} for some constant $c>0$.
 We define $\{V_n(t)\}_{n=1}^\infty$ by
 \begin{align}\label{Vn}
  V_0(t)=v(t); \quad
  V_n(t) = (SV_{n-1})(t), \quad n=1,2,\cdots,
 \end{align}
 where 
 \begin{equation*}
  (Sx)(t):= \int_0^1 G(t,s,\alpha) \overline{h}(s) x(s) ds.
 \end{equation*}
 Then $V_n\in C[0,1]$ and $V_{n-1}(t) \ge V_n(t)\ge 0$ on $[0,1]$ for $n\ge 1$.
 By (iii) of Lemma \ref{propofG}, recalling $||\phi_1||_\infty=1$, we get
 \begin{align*}
  \phi_1(t) = \lambda_1(\alpha) (T\phi_1)(t)
  \le \frac{\lambda_1(\alpha)}{\Gamma(\alpha)} t^{\alpha-1} (1-t) 
        \int_0^1 (1-s)^{\alpha-2} h(s) ds,
 \end{align*}
 where 
 \begin{equation*}
  (Tx)(t):= \int_0^1 G(t,s,\alpha) h(s) x(s) ds.
 \end{equation*}
 Thus, by \eqref{v>ct(1-t)}, there exists $\delta>0$ such that
 $V_0(t) \ge \delta \phi_1(t)$ for $t \in[0,1]$.
 We conclude that
 \begin{align*}
  V_n(t) \ge \delta \phi_1(t), \quad t \in[0,1], \quad n=0,1,2,\cdots,
 \end{align*}
 by induction.
 Indeed, if $V_n(t) \ge \delta \phi_1(t)$ for $t\in[0,1]$, then
 \begin{align*}
  V_{n+1}(t) & = \int_0^1 G(t,s,\alpha) \overline{h}(s) V_n(s) ds \\
   & \ge \lambda_1(\alpha) \int_0^1 G(t,s,\alpha) h(s) \delta \phi_1(s) ds
   = \delta \phi_1(t), \quad t \in [0,1].
 \end{align*}
 Lebesgue's dominated convergence theorem shows that 
 $V(t):=\lim_{n\to\infty}V_n(t)$ satisfies
 $\delta\phi_1(t)\le V(t) \le V_0(t)$ for $t\in[0,1]$ and
 \begin{align*}
  V(t) = (SV)(t) 
  = \int_0^1 G(t,s,\alpha) \overline{h}(s) V(s) ds, \quad t \in [0,1].
 \end{align*}
 Since $(SV)(t)$ is continuous with respect to $t\in[0,1]$, 
 we conclude that $V\in C[0,1]$.
 By Lemma \ref{compact} and (i) of Proposition \ref{KRThm}, 
 we find that $r(S)=1$.
 We set $\widetilde{T}=\lambda_1(\alpha)T$.
 Then 
 \begin{equation*}
  r(\widetilde{T})=r(\lambda_1(\alpha)T)=\lambda_1(\alpha)r(T)=1=r(S).
 \end{equation*}
 Since
 \begin{align*}
  (S-\widetilde{T})(x)(t) 
  & = \int_0^1 G(t,s,\alpha) (\overline{h}(s)-\lambda_1(\alpha)h(s)) x(s) ds,
 \end{align*}
 from (iii) of Lemma \ref{propofG}, for each $x\in K_e \setminus\{0\}$,
 there exists a constant $c_x>0$ such that
 \begin{align*}
  (S-\widetilde{T})(x)(t) \ge c_x t^{\alpha-1}(1-t), \quad t\in[0,1].
 \end{align*}
 This means that $(S-\widetilde{T})(x)\in\mbox{int}(K_e)$ for 
 $x\in K_e \setminus\{0\}$.
 From (iii) of Proposition \ref{KRThm}, it follows that $r(S)>r(\widetilde{T})$.
 This is a contradiction.
\end{proof} 

\begin{lemma}\label{v<lamintGhvhasnosol}
 Let $\alpha\in(1,2]$.
 Assume that \textup{(H)} holds and that $\underline{h}\in C[0,1]$, 
 $\underline{h}(t)\le \lambda_1(\alpha)h(t)$ and 
 $\underline{h}(t)\not\equiv \lambda_1(\alpha)h(t)$ for $t\in[0,1]$.
 Then there is no bounded measurable function $v(t)$ on $[0,1]$ satisfying
 \begin{align*}
  v(t) \le \int_0^1 G(t,s,\alpha) \underline{h}(s) v(s) ds, \quad t \in [0,1]
 \end{align*}
 and \eqref{v>ct(1-t)} for some constant $c>0$.
\end{lemma}

\begin{proof} 
 We suppose that there exists such a function $v(t)$.
 By (iii) of Lemma \ref{propofG}, we get 
 \begin{align*}
  v(t) \le t^{\alpha-1}(1-t) \frac{1}{\Gamma(\alpha)} 
       \int_0^1 (1-s)^{\alpha-2} \underline{h}(s) v(s)ds, \quad t \in [0,1].
 \end{align*}
 By Lemma \ref{lemv>}, there exists a constant $K>0$ such that
 $v(t) \le K \phi_1(t)$ for $t \in [0,1]$.
 We define $\{V_n(t)\}_{n=1}^\infty$ by \eqref{Vn}, where 
 \begin{equation*}
  (Sx)(t):= \int_0^1 G(t,s,\alpha) \underline{h}(s) x(s) ds.
 \end{equation*}
 Then $V_n\in C[0,1]$ and $v(t) \le V_{n-1}(t) \le V_n(t)$ on $[0,1]$ 
 for $n\ge 1$.
 By the same argument as in the proof of Lemma \ref{v>lamintGhvhasnosol},
 we find that $V_n(t)\le K\phi_1(t)$ for $t\in[0,1]$.
 Applying Lebesgue's dominated convergence theorem, we conclude that 
 $V(t):=\lim_{n\to\infty}V_n(t)$ satisfies 
 \begin{align*}
  V(t) = (SV)(t) 
  = \int_0^1 G(t,s,\alpha) \underline{h}(s) V(s) ds, \quad t \in [0,1]
 \end{align*}
 and $v(t)\le V(t) \le K\phi_1(t)$ for $t\in[0,1]$ and that $V\in C[0,1]$.
 Recalling  \eqref{v>ct(1-t)}, we find that 
 $V(t)\ge ct^{\alpha-1}(1-t)$ for $t\in(0,1)$. 
 We set $\widetilde{T}=\lambda_1(\alpha)T$.
 By the same way as in the proof of Lemma \ref{v>lamintGhvhasnosol},
 using Lemma \ref{compact} and (i) of Proposition \ref{KRThm}, we see that 
 $r(S)=r(\widetilde{T})$.
 On the other hand, as in the proof of Lemma \ref{v>lamintGhvhasnosol},
 we conclude that $r(\widetilde{T})>r(S)$.
 This is a contradiction.
\end{proof} 

\begin{proof}[Proof of Theorem \ref{nonexistence}]
 Assume that problem \eqref{P} has a positive solution $u$.
 Recalling Lemma \ref{lemv>}, we note  
 $u(t)\ge (\alpha-1)||u||_{C_{\alpha-2}}t^{\alpha-1}(1-t)$ for $t\in[0,1]$.

 We define the function $\mu(s)$ by
 \begin{align*}
  \mu(s) = \left\{
   \begin{array}{ll}
    s/2, & 0 \le s \le 1, \\[2ex]
    1/2, & s>1.
   \end{array}
  \right.
 \end{align*}
 Then $0\le\mu(s)\le1/2$ for $s \ge 0$ and $\mu(s)/s$ is continuous on 
 $[0,\infty)$.
 We set
 \begin{align*}
  g(s)=\mu(s)f(s) + (1-\mu(s))\lambda_1(\alpha)s, \quad s \ge 0.
 \end{align*} 
 Then $g(s)/s$ is continuous on $[0,\infty)$. 

 We suppose $f(s)/s> \lambda_1(\alpha)$ for $s>0$.
 Then $f(s)/s>g(s)/s>\lambda_1(\alpha)$ for $s>0$.
 We set $\overline{h}(t)=h(t)g(u(t))/u(t)$ for $t\in[0,1]$.
 We see that $\overline{h}\in C[0,1]$, 
 $\overline{h}(t)\ge \lambda_1(\alpha)h(t)$ and 
 $\overline{h}(t)\not\equiv \lambda_1(\alpha)h(t)$ for $t\in[0,1]$.
 Moreover, $u$ satisfies
 \begin{align*}
  u(t) = \int_0^1 G(t,s,\alpha) h(s) f(u(s)) ds
       \ge \int_0^1 G(t,s,\alpha) \overline{h}(s) u(s) ds
 \end{align*}
 for $t\in[0,1]$.
 This contradicts Lemma \ref{v>lamintGhvhasnosol}.

 Similarly, for the case $f(s)/s<\lambda_1(\alpha)$,
 by setting $\underline{h}(t)=h(t)g(u(t))/u(t)$,
 we can lead to a contradiction by Lemma \ref{v<lamintGhvhasnosol}.
\end{proof} 


\section{Existence of positive solutions for the sublinear case} \label{sec:6}

\setcounter{section}{6} \setcounter{equation}{0} 

In this section, we establish the existence of positive solutions for \eqref{P} in the sublinear case, using the Schauder fixed-point theorem.  We define the operators $T$: $C[0,1]\to C[0,1]$ by
\begin{align*}
 (Tu)(t) := \int_0^1 G(t,s,\alpha) h(s) f(|u(s)|) ds,
\end{align*}
for $t\in[0,1]$.
Let $\phi_1$ be an eigenfunction of \eqref{EVP} 
corresponding to $\lambda=\lambda_1(\alpha),$ satisfying $\phi_1(t)>0$ 
for $t\in(0,1)$ and $||\phi_1||_\infty=1$.

\begin{lemma}\label{delta}
 Let $\alpha\in(1,2]$.
 Assume that \textup{(H)} and \textup{(F)} hold and
 that there exists $s_0>0$ such that $f(s)/s>\lambda_1(\alpha)$ for 
 $s\in(0,s_0)$.
 Then there exists $\delta\in(0,1)$ such that
 \begin{align*}
  (T (\delta \phi_1))(t) \ge \delta \phi_1(t), \quad t \in [0,1].
 \end{align*}
\end{lemma}

\begin{proof} 
 Since $f(s)/s>\lambda_1(\alpha)$ for $s\in(0,s_0)$, 
 there exists $\delta\in(0,1)$ such that 
 $f(s) \ge \lambda_1(\alpha)s$ for $s\in[0,\delta]$.
 Therefore, we get
 \begin{equation*}
  (T(\delta\phi_1))(t)
    \ge \int_0^1 G(t,s,\alpha) h(s) \lambda_1(\alpha) \delta\phi_1(s) ds \\
    = \delta \phi_1(t), \quad t\in[0,1].
 \end{equation*}
\end{proof} 

\begin{lemma}\label{M}
 Let $\alpha\in(1,2]$.
 Assume that \textup{(H)} and \textup{(F)} hold and 
 $\displaystyle \limsup_{s\to\infty}f(s)/s<\lambda_1(\alpha)$.
 Then there exists $M>1$ such that
 \begin{align*}
  (T (M \phi_1))(t) \le M \phi_1(t), \quad t \in [0,1].
 \end{align*}
\end{lemma}

\begin{proof} 
 Assume to the contrary that there exist $\{M_n\}_{n=1}^\infty \subset(1,\infty)$ and 
 $\{t_n\}_{n=1}^\infty \subset[0,1]$ such that $M_n\to\infty$ and 
 \begin{equation}
  (T(M_n \phi_1))(t_n)>M_n \phi_1(t_n), \quad n=1,2,\cdots.  
   \label{TMnphi1>}
 \end{equation}
 We note that $t_n\in(0,1)$, that is, $\phi_1(t_n)>0$.
 Indeed, if $t_n=0$ or $1$, then $(T(M_n \phi_1))(t_n)=0$ and $\phi_1(t_n)=0$,
 which is a contradiction.
 We note that the condition $\limsup_{s\to\infty}f(s)/s<\lambda_1(\alpha)$ implies that there exist $\lambda>0$ and $K>0$ such that $0<\lambda<\lambda_1(\alpha)$ and
 \begin{equation*}
  \frac{f(s)}{s} \le \lambda, \quad s\in[K,\infty).
 \end{equation*}
 We set $\overline{f}=\max_{s\in[0,K]}f(s)$.
 From (iii) of Lemma \ref{propofG} and Lemma \ref{lemv>}, it follows that
 \begin{equation*}
  G(t,s,\alpha) 
   \le \frac{1}{\Gamma(\alpha)} t^{\alpha-1} (1-t)(1-s)^{\alpha-2} \\
   \le \frac{\phi_1(t)}{(\alpha-1)\Gamma(\alpha)||\phi_1||_{C_{2-\alpha}}}
 \end{equation*}
 for $t$, $s \in [0,1]$.
 Hence we observe that
 \begin{align*}
  (T(M_n \phi_1))(t_n)
   = & \ \int_{\{s\in[0,1]: 0\le M_n\phi_1(s) < K \}} 
       G(t_n,s,\alpha)h(s)f(M_n\phi_1(s)) ds \\
     & \ + \int_{\{s\in[0,1]: M_n\phi_1(s) \ge K \}} 
       G(t_n,s,\alpha)h(s)f(M_n\phi_1(s)) ds \\
   \le & \ \frac{\overline{f} \phi_1(t_n)}
              {(\alpha-1)\Gamma(\alpha)||\phi_1||_{C_{2-\alpha}}}
        \int_{\{s\in[0,1]: 0\le M_n\phi_1(s) < K \}} h(s) ds \\
     & \ + \lambda M_n \int_{\{s\in[0,1]: M_n\phi_1(s) \ge K \}} 
       G(t_n,s,\alpha)h(s) \phi_1(s) ds \\
   \le & \ \frac{\overline{f} \phi_1(t_n)}
              {(\alpha-1)\Gamma(\alpha)||\phi_1||_{C_{2-\alpha}}}
        \int_0^1 h(s) ds \\
     & + \lambda M_n \int_0^1 G(t_n,s,\alpha)h(s) \phi_1(s) ds \\
   = & \ C \phi_1(t_n) + \frac{\lambda}{\lambda_1(\alpha)} M_n \phi_1(t_n)
 \end{align*}
 for some constant $C>0$.
 Using \eqref{TMnphi1>}, we get
 \begin{equation*}
  M_n \phi_1(t_n) 
  < C \phi_1(t_n) + \frac{\lambda}{\lambda_1(\alpha)} M_n \phi_1(t_n), 
 \end{equation*}
 that is,
 \begin{equation*}
  (\lambda_1(\alpha)-\lambda) M_n < C \lambda_1(\alpha).
 \end{equation*}
 This contradicts the fact that $M_n\to\infty$ and $\lambda<\lambda_1(\alpha)$. 
\end{proof} 

\begin{proof}[Proof of Theorem \ref{existence} for the case \eqref{f0>lam1>foo}]
 We suppose that \eqref{f0>lam1>foo} holds.
 Let $\delta\in(0,1)$ and $M>1$ be constants obtained in Lemmas \ref{delta} 
 and \ref{M}, respectively.
 We define the set $U$ by
 \begin{equation*}
  U = \{ u \in C[0,1] : 
   \delta \phi_1(t) \le u(t) \le M \phi_1(t) \ \mbox{for} \ t \in [0,1] \}.
 \end{equation*}
 Then $U$ is a nonempty bounded convex closed subset of $C[0,1]$.
 Lemmas \ref{delta} and \ref{M} mean that $T$ maps $U$ into itself.
 By the standard argument, we conclude that $T$ is compact and continuous.
 The Schauder fixed-point theorem ensures that there exists $u \in U$ 
 such that $Tu=u$.
 This means that \eqref{P} has a positive solution.
\end{proof} 

\section{Existence of positive solutions for the superlinear case} \label{sec:7}

\setcounter{section}{7} \setcounter{equation}{0} 

In this section, we establish the existence of positive solutions for \eqref{P} in the superlinear case using the Leray-Schauder degree theorem. 
 We define the operator $T$: $C[0,1]\times[0,\infty)\to C[0,1]$ by
\begin{equation*}
 T(u,\eta)(t):= \int_0^1 G(t,s,\alpha) h(s) (f(|u(s)|)+\eta) ds,
  \quad t \in [0,1].
\end{equation*}

\begin{lemma}\label{T(u,eta)not=u}
 Let $\alpha\in(1,2]$.
 Assume that \textup{(H)} and \textup{(F)} hold and
 $\liminf_{s\to\infty}f(s)/s>\lambda_1(\alpha)$.
 Then there exists $\eta>0$ such that
 the equation $T(u,\eta)=u$ has no solution $u\in C[0,1]$.
\end{lemma}

\begin{proof} 
 Assume that, for each $\eta>0$, the equation $T(u,\eta)=u$ has a solution 
 $u\in C[0,1]$.
 Then, for each positive integer $n$, there $\{u_n\}_{n=1}^\infty \subset C[0,1]$ such that 
 $T(u_n,n)=u_n$.
 Lemma \ref{propofG} implies that for $t \in [0,1],$
 \begin{align*}
  u_n(t) \ge \int_0^1 G(t,s,\alpha) h(s) n ds
   \ge \frac{\alpha-1}{\Gamma(\alpha)} t^{\alpha-1}(1-t) n 
    \int_0^1 s(1-s)^{\alpha-1} h(s) ds.
 \end{align*}
 Hence, we see that $u_n(t)\to\infty$ for each fixed $t\in(0,1)$ as $n \to \infty$, so that
 \begin{align*}
  \liminf_{n\to\infty} \frac{f(u_n(t))}{u_n(t)} 
  = \liminf_{s\to\infty} \frac{f(s)}{s} 
 \end{align*}
 for each fixed $t\in(0,1)$.
 We take $\lambda>\lambda_1(\alpha)$ such that
 \begin{align*}
  \liminf_{n\to\infty} \frac{f(u_n(t))}{u_n(t)} \ge \lambda
 \end{align*}
 for each fixed $t\in(0,1)$.
 We set $v_n(t):=u_n(t)/||u_n||_{C_{2-\alpha}}$.
 Then $||v_n||_{C_{2-\alpha}}=1$.
 Lemma \ref{lemv>} implies that 
 $v_n(t)\ge(\alpha-1)t^{\alpha-1}(1-t)$ for $t\in[0,1]$.
 Moreover, $v_n$ satisfies 
 \begin{align*}
  v_n(t) \ge \int_0^1 G(t,s,\alpha) h(s) \frac{f(u_n(s))}{u_n(s)} v_n(s) ds, 
  \quad t\in[0,1].
 \end{align*}
 Using Fatou's lemma, we conclude that $v(t):=\liminf_{n\to\infty}v_n(t)$
 is bounded and measurable on $[0,1]$ (so integrable) and
 \begin{align*}
  v(t) & \ge \liminf_{n\to\infty} \int_0^1 G(t,s,\alpha) h(s) 
           \frac{f(u_n(s))}{u_n(s)} v_n(s) ds \\
  & \ge \int_0^1 G(t,s,\alpha) h(s) 
        \liminf_{n\to\infty} \left(\frac{f(u_n(s))}{u_n(s)} v_n(s) \right) ds \\
  & \ge \int_0^1 G(t,s,\alpha) h(s) 
        \liminf_{n\to\infty} \left( \frac{f(u_n(s))}{u_n(s)} \right)
        \liminf_{n\to\infty} v_n(s) ds \\
  & \ge \lambda \int_0^1 G(t,s,\alpha) h(s) v(s) ds,
  \quad t \in [0,1].
 \end{align*}
 Therefore, Lemma \ref{v>lamintGhvhasnosol} leads us to a contradiction, 
 since $v(t)\ge(\alpha-1)t^{\alpha-1}(1-t)$ for $t\in[0,1]$.
\end{proof} 

\begin{lemma}\label{T(u,eta)not=uM}
 Let $\alpha\in(1,2]$ and $\eta>0$.
 Assume that \textup{(H)} and \textup{(F)} hold and
 $\liminf_{s\to\infty}f(s)/s>\lambda_1(\alpha)$.
 Then there exists $M>1$ such that, for each $\mu\in[0,1]$, the equation
 $T(u,\mu\eta)=u$ has no solution $u\in C[0,1]$ with 
 $||u||_{C_{2-\alpha}}\ge M$.
\end{lemma}

\begin{proof} 
 Assume to the contrary that, for each integer $n>1$, 
 there exists $\mu_n\in[0,1]$ and $u_n\in C[0,1]$ such that 
 $T(u_n,\mu_n\eta)=u_n$ and $||u_n||_{C_{2-\alpha}}\ge n$.
 Lemma \ref{lemv>} implies that $u_n(t)\to\infty$ for each fixed $t\in(0,1)$.
 By the same argument as in the proof of Lemma \ref{T(u,eta)not=u}, 
 we can lead to a contradiction.
\end{proof} 

We define the operator $S$: $C[0,1]\times[0,1]\to C[0,1]$ by
\begin{equation*}
 S(u,\nu)(t):= \nu \int_0^1 G(t,s,\alpha) h(s) f(|u(s)|) ds,
  \quad t \in [0,1].
\end{equation*}

\begin{lemma}\label{S(u,nu)not=u}
 Let $\alpha\in(1,2]$.
 Assume that \textup{(H)} and \textup{(F)} hold and
 that there exists $s_0>0$ such that $f(s)/s<\lambda_1(\alpha)$ for
 $s\in(0,s_0)$.
 Then there exists $\delta\in(0,1)$ such that there is no 
 $(u,\nu)\in(\overline{B(0,\delta)}\setminus\{0\})\times[0,1]$ satisfying 
 $S(u,\nu)=u$, where $B(0,\delta) = \{ u \in C[0,1] : ||u||_\infty < \delta \}$.
\end{lemma}

\begin{proof} 
 Let $\delta\in(0,1)$ satisfy $\delta<s_0$.
 We suppose that there exists 
 $(u,\nu)\in(\overline{B(0,\delta)}\setminus\{0\})\times[0,1]$ such that 
 $S(u,\nu)=u$.
 Since
 \begin{align*}
  u(t) = \nu \int_0^1 G(t,s,\alpha) h(s) f(|u(s)|) ds,  \quad t\in[0,1],
 \end{align*}
 we see that $u(t)\ge 0$ for $t\in[0,1]$ and $\nu\ne0$.
 Moreover, there exists $(t_1,t_2)\subset[0,1]$ such that $t_1<t_2$ and 
 $u(t)>0$ for $t\in(t_1,t_2)$, and hence $||u||_{C_{2-\alpha}}>0$.
 Thus Lemma \ref{lemv>} implies that
 $u(t) \ge ct^{\alpha-1}(1-t)$ for $t \in [0,1]$.
 Now we define a function $g$ by
 \begin{equation*}
  g(s) = \frac{s}{2\delta} f(s) 
     + \left( 1- \frac{s}{2\delta} \right) \lambda_1(\alpha) s,
  \quad s\in [0,\delta].
 \end{equation*}
 Then $g\in C[0,\delta]$ and $f(s)<g(s)<\lambda_1(\alpha)s$ for 
 $s\in(0,\delta]$.
 We set $\underline{h}(t)=h(t)\frac{g(u(t))}{u(t)}$.
 Then $\underline{h}\in C[0,1]$, 
 $\underline{h}(t)\le \lambda_1(\alpha)h(t)$ and 
 $\underline{h}(t)\not\equiv \lambda_1(\alpha)h(t)$ for $t\in[0,1]$, and
 \begin{align*}
  u(t) \le \int_0^1 G(t,s,\alpha) \underline{h}(s) u(s) ds, \quad t\in[0,1].
 \end{align*}
 This contradicts Lemma \ref{v<lamintGhvhasnosol}.
\end{proof} 

\begin{proof}[Proof of Theorem \ref{existence} for the case \eqref{f0<lam1<foo}]
 Assume that \eqref{f0<lam1<foo} holds.
 Then we can take constants $\eta>0$, $M>1$ and $\delta\in(0,1)$ as in Lemmas 
 \ref{T(u,eta)not=u}, \ref{T(u,eta)not=uM} and \ref{S(u,nu)not=u}, respectively.
 We note that $T$ and $S$ are compact and continuous. From Lemma \ref{T(u,eta)not=uM}, it follows that the equation 
 $u-T(u,\mu\eta)=0$ has no solution $(u,\mu)\in \partial B_M\times[0,1]$.
 Thus, $\mbox{deg}_{\rm LS}\,(I-T(\,\cdot\,,\mu\eta),B_M,0)$ 
 is well-defined for $0\le\mu\le1$, where $I$ is an identity, that is, $Iv=v$ for each $v\in C[0,1]$, 
 $\mbox{deg}_{\rm LS}\,(I-T,D,0)$ is the Leray-Schauder degree of $I-T$ in $D$
 over $0$, and $B_M=\{ u\in C[0,1] : ||u||_{C_{2-\alpha}}<M\}$.
 Note that Lemma \ref{T(u,eta)not=u} implies that
 \begin{equation*}
  \mbox{deg}_{\rm LS}\,(I-T(\,\cdot\,,\eta),B_M,0)=0.
 \end{equation*}
 By the photocopy in invariance property of the Leray-Schauder degree,
 we get
 \begin{equation*}
  \mbox{deg}_{\rm LS}\,(I-T(\,\cdot\,,0),B_M,0)
 =\mbox{deg}_{\rm LS}\,(I-T(\,\cdot\,,\eta),B_M,0)=0.
 \end{equation*}
 Let $B(0,\delta)$ be a set as in Lemma \ref{S(u,nu)not=u}.
 We set $B_\delta=\{ u\in C[0,1] : ||u||_{C_{2-\alpha}}<\delta \}$.
 Since $||u||_{C_{2-\alpha}}\le||u||_\infty$, we find $B_\delta\subset B(0,\delta)$.
 Then $\mbox{deg}_{\rm LS}\,(I-S(\,\cdot\,,\nu),B_\delta,0)$ 
 is well-defined for $0\le\nu\le1$, because of Lemma \ref{S(u,nu)not=u}.
 By the homotopy invariance property, we conclude that
 \begin{align*}
      \mbox{deg}_{\rm LS}\,(I-S(\,\cdot\,,1),B_\delta,0) 
    = \mbox{deg}_{\rm LS}\,(I-S(\,\cdot\,,0),B_\delta,0)
    = \mbox{deg}_{\rm LS}\,(I,B_\delta,0)
    = 1.
 \end{align*}
 We note that $T(\,\cdot\,,0)=S(\,\cdot\,,1)$.
 From the excision property of the Leray-Schauder degree, it follows that
 \begin{align*}
  0 & = \mbox{deg}_{\rm LS}\,(I-T(\,\cdot\,,0),B_M,0) \\
    & = \mbox{deg}_{\rm LS}\,(I-T(\,\cdot\,,0),B_\delta,0)
    + \mbox{deg}_{\rm LS}\,(I-T(\,\cdot\,,0),B_M\setminus\overline{B_\delta},0) \\
    & = 1 + \mbox{deg}_{\rm LS}\,(I-T(\,\cdot\,,0),B_M\setminus\overline{B_\delta},0),
 \end{align*}
 which means 
 \begin{equation*}
  \mbox{deg}_{\rm LS}\,(I-T(\,\cdot\,,0),B_M\setminus\overline{B_\delta},0)
  = -1.
 \end{equation*} 
 Consequently, the equation $u-T(u,0)=0$ has a solution
 $u\in B_M\setminus\overline{B_\delta}$, that is,
 $u$ is a positive solution of \eqref{P}.
\end{proof} 

\section{Uniqueness of positive solutions for the sublinear case} \label{sec:8}

\setcounter{section}{8} \setcounter{equation}{0} 

In this section, we prove Theorem \ref{uniquenesssub}, which establishes the uniqueness of positive solutions for \eqref{P} in the sublinear case.

\begin{proof}[Proof of Theorem \ref{uniquenesssub}]
By recalling Remark \ref{remoff(s)/s}, Theorem \ref{nonexistence} guarantees that problem \eqref{P} has at least one positive solution.

 We suppose that problem \eqref{P} has two positive solutions $u$ and $v$.
 From (iii) of Lemma \ref{propofG}, it follows that 
 \begin{align}
  u(t) \le t^{\alpha-1}(1-t) \frac{1}{\Gamma(\alpha)} ||f(u)||_\infty
       \int_0^1 (1-s)^{\alpha-2} h(s) ds, \quad t \in [0,1].
   \label{u(t)<Ce(t)}
 \end{align}
 Hence, by Lemma \ref{lemv>}, there exists a constant $\beta_1\ge 1$ such that
 $u(t) \le \beta_1 v(t)$ for $t \in [0,1]$.
 We set
 \begin{align*}
  \beta_0 = \inf \{ \beta\ge1 : 
  u(t) \le \beta v(t) \ \textup{for} \ t \in [0,1] \}.
 \end{align*}
 Then $1\le\beta_0\le\beta_1$.
 We will prove $\beta_0=1$.
 Assume $\beta_0>1$.
 We define $w(t)=\beta_0 v(t)$.
 We have $u(t)\le w(t)$ and $v(t)<w(t)$ for $t\in(0,1)$.
 Since
 \begin{align*}
  f(w(t)) = \frac{f(w(t))}{w(t)} w(t) < \frac{f(v(t))}{v(t)} w(t)
  = \beta_0f(v(t)), \quad t\in(0,1), 
 \end{align*}
 we get
 \begin{align}\label{w>intGhf(w)}
  w(t) = \beta_0 v(t) & = \beta_0 \int_0^1 G(t,s,\alpha) h(s) f(v(s))ds \\
  & > \int_0^1 G(t,s,\alpha) h(s) f(w(s))ds, \quad t \in (0,1).  \nonumber
 \end{align}
 By Lemma \ref{propofG}, we see that 
 \begin{align*}
  w(t)-u(t) & \ge \int_0^1 G(t,s,\alpha) h(s) f(w(s))ds
                 - \int_0^1 G(t,s,\alpha) h(s) f(u(s))ds \\
   & = \int_0^1 G(t,s,\alpha) h(s) (f(w(s))-f(u(s))) ds
     \ge K t^{\alpha-1}(1-t)
 \end{align*}
 for $t\in[0,1]$, where
 \begin{align*}
  K = \frac{\alpha-1}{\Gamma(\alpha)}  
      \int_0^1 s(1-s)^{\alpha-1} h(s) (f(w(s))-f(u(s))) ds \ge 0.
 \end{align*}
 Now we assume that $w(t_1)>u(t_1)$ for some $t_1\in(0,1)$.
 Then $w(t)>u(t)$ for $t\in(t_1-\delta,t_1+\delta)$ for some $\delta>0$, 
 which implies $K>0$.
 In the same way of \eqref{u(t)<Ce(t)}, there exists $C>0$ such that
 $v(t)\le Ct^{\alpha-1}(1-t)$ for $t\in[0,1]$, and hence 
 \begin{align*}
  u(t) \le w(t)- Kt^{\alpha-1}(1-t) 
       \le \beta_0 v(t) - \frac{K}{C} v(t)
       = \Bigr( \beta_0 - \frac{K}{C} \Bigl) v(t)
 \end{align*}
 for $t\in[0,1]$.
 This contradicts the definition of $\beta_0$.
 Therefore, $w(t)=u(t)$ for $t\in[0,1]$.
 By \eqref{w>intGhf(w)}, we conclude that
 \begin{align*}
  u(t) = w(t) > \int_0^1 G(t,s,\alpha) h(s) f(w(s))ds
   = \int_0^1 G(t,s,\alpha) h(s) f(u(s))ds = u(t)
 \end{align*}
 for $t \in (0,1)$.
 This is a contradiction.
 We find $\beta_0=1$, that is, $u(t)\le v(t)$ for $t\in[0,1]$.
 Reversing the roles of $u$ and $v$, we get $v(t)\le u(t)$ for $t\in[0,1]$.
 Consequently, $u(t)=v(t)$ for $t\in[0,1]$.
\end{proof} 

\section{Uniqueness of positive solutions for the superlinear case} \label{sec:9}

\setcounter{section}{9} \setcounter{equation}{0} 

In this section, we prove the following result and, by applying it, provide a proof of Theorem \ref{uniquenesssuper}. We recall that a solution $u$
of \eqref{P} is said to be {\it nondegenerate} if the following problem has no nontrivial solution:
\begin{align}\label{L}
 \left\{
  \begin{array}{l}
   D_{0+}^{\alpha} w + h(t)f'(u(t)) w = 0, \quad 0<t<1, \\[1ex]
   w(0)=w(1)=0.
  \end{array}
 \right.
\end{align}

\begin{theorem}\label{atmost}
 Let $\alpha_0\in(1,2]$ and let $k$ be a positive integer.
 Assume that \textup{(H)} and \textup{(F)} hold, $f\in C^1[0,\infty)$,
 either $\liminf_{s\to0^+}f(s)/s>\lambda_1(\alpha_0)$ or \linebreak
 $\limsup_{s\to0^+}f(s)/s<\lambda_1(\alpha_0)$,
 and that one of the following conditions holds\textup{:}
 \begin{enumerate}
  \item $\liminf_{s\to\infty}f(s)/s > \lambda_1(\alpha_0)$ and 
	$\limsup_{s\to\infty}s^{-p}f(s)<\infty$ for some $p>0$ with 
	$p(2-\alpha_0)<1$\textup{;}
  \item $\limsup_{s\to\infty}f(s)/s < \lambda_1(\alpha_0)$.
 \end{enumerate} 
 If problem \eqref{P} with $\alpha=\alpha_0$ has at most $k$ positive solutions 
 and every positive solution of \eqref{P} with $\alpha=\alpha_0$ is 
 nondegenerate, then there exists $\rho>0$ such that, 
 for each $\alpha\in(\alpha_0-\rho,\alpha_0+\rho)\cap(1,2]$, 
 problem \eqref{P} has at most $k$ positive solutions.
\end{theorem}

\begin{lemma}\label{||u||<general1}
 Let $\alpha_*$ and $\alpha^*$ satisfy $1<\alpha_*\le\alpha^*\le 2$. 
 Assume that \textup{(H)} and \textup{(F)} hold and
 \begin{align}
  \liminf_{s\to\infty} \frac{f(s)}{s}
  > \sup_{\alpha\in[\alpha_*,\alpha^*]} \lambda_1(\alpha).
  \label{f(s)/s>suplam1}
 \end{align}
 Then there exists a constant $M>0$ such that
 \begin{align*}
  ||u||_{C_{2-\alpha}} \le M
 \end{align*}
 for every positive solution $u$ of \eqref{P} and 
 $\alpha\in[\alpha_*,\alpha^*]$,
 where the constant $M$ does not depend on $u$ and $\alpha$. 
\end{lemma}

\begin{proof} 
 The proof is similar to Lemma \ref{T(u,eta)not=uM}.
 Assume that, for each integer $n\ge 1$, there exist $u_n\in C[0,1]$ and 
 $\alpha_n\in[\alpha_*,\alpha^*]$ such that $u_n$ is a positive solution of 
 \eqref{P} with $\alpha=\alpha_n$ and $||u_n||_{C_{2-\alpha_n}}>n$.
 Lemma \ref{lemv>} ensures that $u_n(t)\to\infty$ for each fixed $t\in(0,1)$.
 We take $\lambda$ for which 
 \begin{align*}
  \liminf_{s\to\infty} \frac{f(s)}{s} > \lambda
  > \sup_{\alpha\in[\alpha_*,\alpha^*]} \lambda_1(\alpha).
 \end{align*}
 Then
 \begin{align*}
  \liminf_{n\to\infty} \frac{f(u_n(t))}{u_n(t)} 
  = \liminf_{s\to\infty} \frac{f(s)}{s} > \lambda
 \end{align*}
 for each fixed $t\in(0,1)$.
 In the same way of Lemma \ref{T(u,eta)not=u}, we deduce that
 $v_n(t):=u_n(t)/||u_n||_{C_{2-\alpha}}$ satisfies $||v_n||_{C_{2-\alpha}}=1$,
 \begin{align*}
  v_n(t)\ge(\alpha_*-1)t^{\alpha^*-1}(1-t), \quad t\in[0,1]
 \end{align*}
 and
 \begin{align*}
  v_n(t) = \int_0^1 G(t,s,\alpha_n) h(s) \frac{f(u_n(s))}{u_n(s)} v_n(s) ds, 
  \quad t\in[0,1].
 \end{align*}
 Since $\{\alpha_n\}_{n=1}^\infty\subset [\alpha_*,\alpha^*]$, 
 there exists a convergent subsequence of $\{\alpha_n\}_{n=1}^\infty$.
 We denote it by $\{\alpha_n\}$ again and its limit by 
 $\alpha_0\in[\alpha_*,\alpha^*]$.
 By using the same arguments as in the proof of Lemma \ref{T(u,eta)not=u} again,
 Fatou's lemma shows that $v(t):=\liminf_{n\to\infty}v_n(t)$
 is measurable on $[0,1]$ and
 \begin{align*}
  v(t) = \liminf_{n\to\infty} \int_0^1 G(t,s,\alpha_n) h(s) 
           \frac{f(u_n(s))}{u_n(s)} v_n(s) ds
   \ge \lambda \int_0^1 G(t,s,\alpha_0) h(s) v(s) ds
 \end{align*}
 for $t \in [0,1].$
 This contradicts Lemma \ref{v>lamintGhvhasnosol}, because of 
 $v(t)\ge(\alpha_*-1)t^{\alpha^*-1}(1-t)$ for $t\in[0,1]$.
\begin{proof}end 

\begin{lemma}\label{||u||<general2}
 Let $\alpha_*$ and $\alpha^*$ satisfy $1<\alpha_*\le\alpha^*\le 2$. 
 Assume that \textup{(H)} and \textup{(F)} hold and
 \begin{align}
  \limsup_{s\to\infty} \frac{f(s)}{s}
  <\inf_{\alpha\in[\alpha_*,\alpha_*]} \lambda_1(\alpha).
  \label{f(s)/s<inflam1}
 \end{align}
 Then there exists a constant $M>0$ such that
 \begin{align*}
  ||u||_{C_{2-\alpha}} \le M
 \end{align*}
 for every positive solution $u$ of \eqref{P} and 
 $\alpha\in[\alpha_*,\alpha^*]$,
 where the constant $M$ does not depend on $u$ and $\alpha$. 
\end{lemma}

\begin{proof} 
 We suppose that, for each integer $n\ge 1$, there exists $u_n\in C[0,1]$ and 
 $\alpha_n\in[\alpha_*,\alpha^*]$ such that $u_n$ is a positive solution of 
 \eqref{P} with $\alpha=\alpha_n$ and $||u_n||_{C_{2-\alpha_n}}>n$.
 Let $\lambda$ satisfy
 \begin{align*}
    \limsup_{s\to\infty} \frac{f(s)}{s} < \lambda
  < \inf_{\alpha\in[\alpha_*,\alpha_*]} \lambda_1(\alpha)
 \end{align*}
 From Lemma \ref{lemv>}, it follows that $u_n(t)\to\infty$ for each fixed 
 $t\in(0,1)$, which implies
 \begin{align*}
    \limsup_{n\to\infty} \frac{f(u_n(t))}{u_n(t)} 
  = \limsup_{s\to\infty} \frac{f(s)}{s} 
  < \lambda
 \end{align*}
 for each fixed $t\in(0,1)$.
 We set $v_n(t):=u_n(t)/||u_n||_{C_{2-\alpha}}$.
 Since $||v_n||_{C_{2-\alpha}}=1$, it holds that
 \begin{align*}
  v_n(t) \le t^{\alpha_n-2} \le t^{\alpha_*-2}, \quad t \in (0,1]
 \end{align*}
 and
 \begin{align*}
  v_n(t)\ge(\alpha_*-1)t^{\alpha^*-1}(1-t), \quad t\in[0,1],
 \end{align*}
 because of Lemma \ref{lemv>}.
 We note here that $t^{\alpha_*-2}$ is integrable on $[0,1]$.
 We take a convergent subsequence of $\{\alpha_n\}_{n=1}^\infty$.
 We denote it by $\{\alpha_n\}$ again and its limit by 
 $\alpha_0\in[\alpha_*.\alpha^*]$.
 Applying reverse Fatou's lemma, we find that 
 $v(t):=\limsup_{n\to\infty}v_n(t)$ is measurable on $[0,1]$ and
 \begin{align*}
  v(t) & = \limsup_{n\to\infty} \int_0^1 G(t,s,\alpha_n) h(s) 
           \frac{f(u_n(s))}{u_n(s)} v_n(s) ds \\
  & \le \int_0^1 h(s) 
        \limsup_{n\to\infty} G(t,s,\alpha_n) \frac{f(u_n(s))}{u_n(s)} v_n(s) ds \\
  & \le \int_0^1 h(s) \limsup_{n\to\infty} G(t,s,\alpha_n)  
        \limsup_{n\to\infty} \frac{f(u_n(s))}{u_n(s)} 
        \limsup_{n\to\infty} v_n(s) ds \\
  & \le \lambda \int_0^1 G(t,s,\alpha_0) h(s) v(s) ds, \quad t \in [0,1].
 \end{align*}
 This contradicts Lemma \ref{v<lamintGhvhasnosol}.
\end{proof} 

\begin{lemma}\label{un->u}
 Let $\alpha_0\in(1,2]$ and $\{\alpha_n\}_{n=1}^\infty$ satisfy 
 $\alpha_n\in(1,2]$ for $n\ge1$ and $\alpha_n\to \alpha_0$ as $n\to\infty$.
 Suppose that \textup{(H)} and \textup{(F)} hold and
 either $\displaystyle \liminf_{s\to0^+}f(s)/s>\lambda_1(\alpha_0)$ or
 $\displaystyle \limsup_{s\to0^+}f(s)/s<\lambda_1(\alpha_0)$.
 For each integer $n\ge1$, let $u_n$ be a positive solution of \eqref{P} 
 with $\alpha=\alpha_n$.
 Assume, moreover, that $u_n$ converges to $u_0$ uniformly on $[0,1]$ as 
 $n\to\infty$ for some $u_0 \in C[0,1]$.
 Then $u_0$ is a positive solution of \eqref{P} with $\alpha=\alpha_0$.
\end{lemma}

\begin{proof} 
 First we prove that $u_0$ is a nonnegative solution of \eqref{P} with 
 $\alpha=\alpha_0$.
 Since $u_n(t)>0$ for $t\in(0,1)$, we find that $u_0(t)\ge 0$ for $t\in[0,1]$.
 Since $u_n$ converges to $u_0\in C[0,1]$ uniformly on $[0,1]$, 
 there exists $M>0$ such that $||u_n||_\infty \le M$ for $n\ge1$.
 From (ii) of Lemma \ref{propofG}, it follows that
 \begin{align}\label{G<K}
  0 \le G(t,s,\alpha) \le \frac{1}{\Gamma(\alpha)}(s(1-s))^{\alpha-1}
    \le L, \ \ t,s \in [0,1], \ \alpha\in(1,2]
 \end{align}
 for some constant $L>0$. 
 By letting $n\to\infty$ in
 \begin{align}
  u_n(t) = \int_0^1 G(t,s,\alpha_n) h(s) f(u_n(s)) ds, 
  \quad t\in[0,1].
   \label{u_n=int}
 \end{align}
 Lebesgue's dominated convergence theorem yields
 \begin{align}\label{u_0}
  u_0(t) = \int_0^1 G(t,s,\alpha_0) h(s) f(u_0(s)) ds, 
  \quad t\in[0,1].
 \end{align}
 This means that $u_0$ is a nonnegative solution of \eqref{P} with
 $\alpha=\alpha_0$.

 Next we claim that $u_0(t_0)>0$ for some $t_0\in(0,1)$.
 Assume to the contrary that $u_0(t)=0$ for $t\in[0,1]$.
 We suppose that $\liminf_{s\to0^+}f(s)/s>\lambda_1(\alpha_0)$.
 By Proposition \ref{contioflam1}, there exists $\rho>0$ such that
 \begin{align*}
  \liminf_{s\to0^+}\frac{f(s)}{s} > \max_{\alpha\in I}\lambda_1(\alpha)
 \end{align*}
 where $I=[\alpha_0-\rho,\alpha_0+\rho]\cap(1,2]$.
 Then there exist $\delta>0$ and $\lambda>\max_{\alpha\in I}\lambda_1(\alpha)$ 
 such that $f(s)/s\ge \lambda$ for $s\in(0,\delta]$.
 Since $u_n$ converges to $0$ uniformly on $[0,1]$,
 there exist an integer $N>0$ such that $0<u_N(t)<\delta$ for $t\in(0,1)$
 and $\alpha_N\in I$, and hence 
 \begin{align*}
  \frac{f(u_N(t))}{u_N(t)} \ge \lambda > \lambda_1(\alpha_N), \quad t\in(0,1).
 \end{align*}
 By \eqref{u_n=int}, we have
 \begin{align*}
  u_N(t) \ge \lambda \int_0^1 G(t,s,\alpha_N) h(s) u_N(s) ds, 
  \quad t\in[0,1].
 \end{align*}
 By recalling Lemma \ref{lemv>}, this contradicts 
 Lemma \ref{v>lamintGhvhasnosol}.
 Similarly, in the case where $\limsup_{s\to0^+}f(s)/s<\lambda_1(\alpha_0)$, 
 by using Lemma \ref{v<lamintGhvhasnosol}, we can lead to a contradiction.
 Therefore, $u_0(t_0)>0$ for some $t_0\in(0,1)$ as claimed.

 Since $u_0$ satisfies \eqref{u_0}, Lemma \ref{lemv>} implies that 
 $u_0(t)>0$ for $t\in(0,1)$.
 Consequently, $u_0$ is a positive solution of \eqref{P} with $\alpha=\alpha_0$.
\end{proof} 

\noindent{\it Proof of Theorem \ref{atmost}}
 Assume to the contrary that there exist 
 $\{\alpha_n\}_{n=1}^\infty$ and $\{u_{i,n}\}_{n=1}^\infty$, 
 $i=1$, $2$, $\cdots$, $k+1$ such that $\alpha_n\in(1,2]$ for $n\ge1$, 
 $\alpha_n \to \alpha_0$ as $n \to \infty$, 
 $u_{1,n}$, $u_{2,n}$ $\cdots$, $u_{k+1,n}$ are positive solutions of \eqref{P}
 with $\alpha=\alpha_n$ different from each other.
 By Proposition \ref{contioflam1}, there exist $\alpha_*$, $\alpha^*$ and $N>0$ 
 such that $1<\alpha_*<\alpha_0<\alpha^*\le2$, $\alpha_n\in[\alpha_*,\alpha^*]$
 for $n\ge N$ and either \eqref{f(s)/s>suplam1} or \eqref{f(s)/s<inflam1} holds.
 Moreover, when (i) holds, we take $\alpha_*$ for which $p(2-\alpha_*)<1$. 
 Lemmas \ref{||u||<general1} and \ref{||u||<general2} ensure that
 there exists a constant $M>0$ such that 
 \begin{align*}
  0<u_{i,n}(t) \le Mt^{\alpha_*-2}, \quad t \in (0,1), \ 
  i\in\{1,2,\cdots,k+1\}, \ n \ge N.
 \end{align*}
 and $M$ does not depend on $i$ and $n$.
 We note that $u_{i,n}(0)=u_{i,n}(1)=0$.
 Then $\{u_{i,n}(t)\}_{n=1}^\infty$, $i=1$, $2$, $\cdots$, $k+1$ are bounded
 at every $t\in[0,1]$.
 
 Now we will prove that $\{u_{i,n}(t)\}_{n=1}^\infty$, $i=1$, $2$, $\cdots$, 
 $k+1$ are equicontinuos on $[0,1]$. 
 First we assume that (i) holds.
 Then there exist $s_1>0$ and $c_1>0$ such that $s^{-p}f(s)\le c_1$ for 
 $s\ge s_1$, and hence
 \begin{align*}
  f(s) \le c_1 s^p + c_2, \quad s\ge 0,
 \end{align*}
 where $c_2:=\max_{t\in[0,s_1]}f(t)$.
 Let $\varepsilon>0$ be arbitrary.
 Since $G(t,s,\alpha)$ is uniformly continuous on 
 $[0,1]\times[0,1]\times[\alpha_*,\alpha^*]$, there exists $\delta>0$ such that
 if $t_1$, $t_2\in[0,1]$, $|t_1-t_2|<\delta$, $s\in[0,1]$, 
 $\alpha\in[\alpha_*,\alpha^*]$, then
 \begin{align*}
  |G(t_1,s,\alpha)-G(t_2,s,\alpha)|< \varepsilon.
 \end{align*}
 If $|t_1-t_2|<\delta$, then
 \begin{align*}
  |u_{i,n}(t_1)-u_{i,n}(t_2)|
   & = \left| \int_0^1 (G(t_1,s,\alpha_n)-G(t_2,s,\alpha_n)) h(s) 
              f(u_{i,n}(s)) ds \right| \\
   & \le ||h||_\infty 
   \int_0^1 |G(t_1,s,\alpha_n)-G(t_2,s,\alpha_n)| (c_1(u_{i,n}(s))^p+c_2) ds \\
   & < \varepsilon ||h||_\infty \int_0^1 (c_1 M^p s^{p(\alpha_*-2)}+c_2) ds \\
   & = \varepsilon ||h||_\infty 
       \left( \frac{c_1 M^p}{1-p(2-\alpha_*)} + c_2 \right).
 \end{align*}
 Therefore, $\{u_{i,n}(t)\}_{n=1}^\infty$ is equicontinuos.
 Now we suppose that (ii) holds.
 Then there exists $s_1>0$ such that $f(s)/s<\lambda_1(\alpha_0)$ for  
 $s\ge s_1$, which yields
 \begin{align*}
  f(s) \le \lambda_1(\alpha_0) s + \max_{t\in[0,s_1]}f(t), \quad s\ge 0.
 \end{align*}
 In the same way as (i), we can deduce that
 $\{u_{i,n}(t)\}_{n=1}^\infty$ is equicontinuos. 

 Hence, the Arzel\`{a}-Ascoli theorem implies that 
 $\{u_{i,n}(t)\}_{n=1}^\infty$ contains a subsequences converging 
 uniformly on $[0,1]$.
 For convenience, we use the same notation $\{u_{i,n}\}_{n=1}^\infty$,
 $i=1$, $2$, $\cdots$, $k+1$ for the uniformly convergent subsequences.
 Namely, $\{u_{i,n}\}_{n=1}^\infty$ converges uniformly to $u_i$ 
 for some $u_i\in C[0,1]$.
 Lemma \ref{un->u} implies that $u_i$, $i=1,2,\cdots,k+1$ are positive 
 solutions of \eqref{P} with $\alpha=\alpha_0$.
 Since problem \eqref{P} with $\alpha=\alpha_0$ has at most $k$ positive 
 solutions, we find that $u_l=u_m$ for some $l$, $m\in \{1,2,\cdots,k+1\}$ 
 with $l\ne m$.
 Now we set $U:=u_l=u_m$.
 We note that $U$ is nondegenerate by the assumption. 
 We set
 \[
  w_n(t) := \frac{u_{l,n}(t)-u_{m,n}(t)}{||u_{l,n}-u_{m,n}||_\infty}.
 \]
 Then $||w_n||_\infty=1$ and $w_n(0)=w_n(1)=0$.
 It holds that
 \begin{equation*}
  w_n(t) = \int_0^1 G(t,s,\alpha_n) h(s)
   \frac{f(u_{l,n}(s))-f(u_{m,n}(s))}{||u_{l,n}-u_{m,n}||_\infty} ds,
   \quad t \in [0,1].
 \end{equation*}
 By the mean value theorem, for each $s\in(0,1)$ and $n\ge N$, 
 there exists $z_n(s)$ such that 
 \begin{align*}
  f(u_{l,n}(s))-f(u_{m,n}(s)) = f'(z_n(s))(u_{l,n}(s)-u_{m,n}(s))
 \end{align*}
 and 
 \begin{align*}
  \min\{u_{l,n}(s),u_{m,n}(s)\} \le z_n(s) \le \max\{u_{l,n}(s),u_{m,n}(s)\}. 
 \end{align*}
 Since $u_{l,n}$ and $u_{m,n}$ converge to $U\in C[0,1]$ uniformly on $[0,1]$,
 we find that $\{u_{l,n}\}_{n=1}^\infty$ and $\{u_{m,n}\}_{n=1}^\infty$ 
 are uniformly bounded on $[0,1]$.
 Hence, there exists $C>0$ such that $|f'(z_n(s))| \le C$ 
 for $s\in[0,1]$ and $n\ge N$.
 We observe that
 \begin{align*}
  \left| \frac{f(u_{l,n}(s))-f(u_{m,n}(s))}{||u_{l,n}-v_{m,n}||_\infty} \right|
  & = |f'(z_n(s))| |w_n(s)| \le C, \quad s \in [0,1]. 
 \end{align*}
 By the same argument above, we find that $\{w_n\}_{n=1}^\infty$ is 
 equicontinuos on $[0,1]$. 
 By the Arzel\`{a}-Ascoli theorem, there exists a subsequence of 
 $\{w_n\}_{n=1}^\infty$ that converges to $w$ uniformly on $[0,1]$ for some $w \in C[0,1]$.
 We use the same notation $\{w_n\}$ for the uniformly convergent subsequence.
 Therefore, we conclude that
 \begin{align*}
  \frac{f(u_{l,n}(t))-f(u_{m,n}(t))}{||u_{l,n}-v_{m,n}||_\infty}
  = f'(z_n(s)) w_n(s) \to f'(U(s))w(s)
  \quad \mbox{as} \ n \to \infty
 \end{align*}
 for each fixed $s \in [0,1]$.
 Lebesgue's dominated convergence theorem implies 
 \begin{align*}
  w(t) = \int_0^1 G(t,s,\alpha_0) h(s)f'(U(s)) w(s) ds, \quad t \in [0,1],
 \end{align*}
 which means that $w$ is a nontrivial solution of \eqref{L} with 
 $\alpha=\alpha_0$ and $u=U$, by recalling $||w_n||_\infty=1$.
 On the other hand, since $U$ is a positive solution of \eqref{P} with 
 $\alpha=\alpha_0$, it is nondegenerate, that is,
 problem \eqref{L} with $\alpha=\alpha_0$ and $u=U$ has no nontrivial solution.
 This is a contradiction.
\end{proof} 

Now we prove Theorem \ref{uniquenesssuper}.
We consider the case $\alpha=2$, that is,
\begin{align}\label{P2}
 \left\{
  \begin{array}{l}
   u'' + h(t)f(u) = 0, \quad 0<t<1, \\[1ex]
   u(0)=u(1)=0.
  \end{array}
 \right.
\end{align}
We define the function $\overline{f}$ by
\begin{align*}
 \overline{f}(s) := \left\{ 
  \begin{array}{ll}
   f(s), & s\ge 0, \\[1ex]
   -f(-s), & s<0.
  \end{array}
 \right.
\end{align*}
Recalling Remark \ref{remoff}, we find that
if $f\in C^1[0,\infty)$, $f(s)>0$ for $s>0$ and \eqref{f'>f/s} holds, then
$\overline{f}\in C^1(\mathbb{R})$, $s\overline{f}(s)>0$ and 
$\overline{f}'(s)>\overline{f}(s)/s$ for $s\ne0$.
The following result is obtained in \cite[Theorem 1.1]{T1}.

\begin{proposition}\label{uniqueness2}
 Suppose that $h\in C^1[0,1]$, $h(t)>0$ for $t\in[0,1]$, and
 \eqref{<h'/h<} holds.
 Assume, moreover, that $f\in C^1[0,\infty)$, $f(s)>0$ for $s>0$, and
 \eqref{f'>f/s} holds.
 Then problem \eqref{P2} has at most one positive solution.
\end{proposition}

\begin{lemma}\label{nondegerate2}
 Under the condition of Proposition \ref{uniqueness2},
 if problem \eqref{P2} has a positive solution $u$, then $u$ is nondegenerate.
\end{lemma}

\begin{proof} 
 Let $w$ be a solution of initial value problem
 \begin{align*}
  \left\{
  \begin{array}{l}
   w''+h(t)f'(u(t)) w = 0, \quad t \in (0,1) \\[1ex]
   w(0)=0, \ w'(0)=1.
  \end{array}
  \right.
 \end{align*}
 Then Lemmas 2.2, 2.3 and 2.4 in \cite{T1} imply that
 $w(t)$ has a unique zero in $(0,1)$ and $w(1)<0$.
 This means that $u$ is nondegenerate.
\end{proof} 

\noindent{\it Proof of Theorem \ref{uniquenesssuper}}
 From Remark \ref{remoff}, it follows that \eqref{f0<lam1<foo} with $\alpha=2$ 
 holds.
 By Lemma \ref{contioflam1} and Theorem \ref{existence}, there exists 
 $\rho_1\in(1,2)$ such that problem \eqref{P} has a positive solution 
 for each fixed $\alpha\in(2-\rho_1,2]$.
 By Proposition \ref{uniqueness2} and Lemma \ref{nondegerate2},
 we conclude that problem \eqref{P} with $\alpha=2$ has a unique positive 
 solution and it is nondegenerate.
 Consequently, by Theorem \ref{atmost} with $\alpha_0=2$ and $k=1$, 
 there exists $\rho\in(0,\rho_1]$ such that problem \eqref{P} has a unique 
 positive solution for each fixed $\alpha\in(2-\rho,2]$.
\end{proof} 

\section{Existence of multiple positive solutions} \label{sec:10}

\setcounter{section}{10} \setcounter{equation}{0} 

In this section, we provide a proof of Theorem \ref{nonuniquenessexample}. To this end, we prove the following existence result using the results from Section \ref{sec:3}.

\begin{theorem}\label{existencealphanear2}
 Let $\alpha_0\in(1,2]$ and let $k$ be a positive integer.
 Suppose that \textup{(H)} and \textup{(F)} hold and $f\in C^1[0,\infty)$.
 If problem \eqref{P} with $\alpha=\alpha_0$ has at least $k$ nondegenerate 
 positive solution, then there exists $\rho>0$ such that, 
 for each $\alpha\in(\alpha_0-\rho,\alpha_0+\rho)\cap(1,2]$, 
 problem \eqref{P} has at least $k$ positive solutions. 
\end{theorem}

\begin{proof} 
 We take a function $\overline{f}$ for which 
 $\overline{f}\in C^1(\mathbb{R})$, $\overline{f}(s)=f(s)$ for $s\ge0$,
 $f(0)-1\le \overline{f}(s)\le f(0)+1$ for $s\in[-1,0]$, and
 $\overline{f}(s)=f(0)$ for $s\le -1$.
 We define the operator $T$: $X_e\times(1,2]\to X_e$ by 
 \begin{align*}
  T(u,\alpha)(t) = u(t) - \int_0^1 G(t,s,\alpha) h(s) \overline{f}(u(s)) ds, 
  \quad t\in[0,1].
 \end{align*}
 We observe that 
 \begin{align*}
  (D_u T(u,\alpha)w)(t)
   = w(t) - \int_0^1 G(t,s,\alpha) h(s) \overline{f}'(u(s)) w(s) ds,
   \quad t\in[0,1].
 \end{align*}
 Let $u_0$ be a nondegenerate positive solution of \eqref{P} with
 $\alpha=\alpha_0$.
 Then
 \begin{align*}
  (D_u T(u_0,\alpha_0)w)(t)  
  = w(t) - \int_0^1 G(t,s,\alpha_0) h(s) f'(u_0(s)) w(s) ds.
 \end{align*}
 Since $u_0$ is nondegenerate, we find that $D_u T(u_0,\alpha_0)$ is injective.
 We define the operator $S$: $X_e \to X_e$ by
 \begin{align*}
  (Sw)(t):= \int_0^1 G(t,s,\alpha_0) h(s) f'(u_0(s))w(s) ds, 
  \quad t\in[0,1].
 \end{align*}
 Then $D_u T(u_0,\alpha_0)=I-S$, where $I$ is the identity map, 
 that is, $Iw=w$ for every $w\in X_e$.
 Lemma \ref{compact} implies that $S$ is compact.
 From the Fredholm alternative, it follows that $D_u T(u_0,\alpha_0)=I-S$ is 
 homeomorphism.
 By the implicit function theorem, there exist $\rho_0>0$ and a unique 
 continuous mapping 
 $u$: $(\alpha_0-\rho_0,\alpha_0+\rho_0)\cap(1,2] \to X_e$ such that
 \begin{align*}
  T(u(\alpha),\alpha)=0, \quad 
  \alpha \in (\alpha_0-\rho_0,\alpha_0+\rho_0)\cap(1,2]
 \end{align*}
 and $||u(\alpha)-u_0||_e<\rho_0$, $u(\alpha_0)=u_0$.
 Then $u(\alpha)$ satisfies 
 \begin{align*}
  u(\alpha)(t) = \int_0^1 G(t,s,\alpha) h(s) \overline{f}(u(\alpha)(s)) ds, 
  \quad t\in[0,1].
 \end{align*}
 By Lemma \ref{lemv>}, We conclude that
 \begin{align*}
  u(\alpha)(t)  
  & = u(\alpha)(t) - u_0(t) + u_0(t) \\
  & \ge -||u(\alpha)-u_0||_e e(t) + (\alpha_0-1) ||u_0||_{C_{2-\alpha_0}} e(t) \\
  & = ((\alpha_0-1)||u_0||_{C_{2-\alpha_0}}-||u(\alpha)-u_0||_e) e(t),
  \quad t \in [0,1],
 \end{align*}
 where $e(t)=t^{\alpha-1}(1-t)$.
 Since $||u(\alpha)-u_0||_e\to0$ as $\alpha\to\alpha_0$,
 there exists $\rho_1 \in (0,\rho_0]$ such that
 if $\alpha\in(\alpha_0-\rho_1,\alpha_0+\rho_1)\cap(1,2]$, 
 then $u(\alpha)(t)>0$ for $t\in(0,1)$,
 which means that $u(\alpha)$ is a positive solution of \eqref{P}.
 Consequently, there exists $\rho>0$ such that, for each 
 $\alpha\in(\alpha_0-\rho,\alpha_0+\rho)\cap(1,2]$, 
 problem \eqref{P} has at least $k$ positive solutions.
\end{proof} 

The proof of the following result will be given in the next section.

\begin{lemma}\label{nonuniqueness2}
 Let $l>1$ and $p>1$ satisfy $(p-1)l\ge 4$.
 Then there exists $\delta_0>0$ such that 
 \begin{align}
  \left\{
  \begin{array}{l}
   u'' + \left| t -\dfrac{1}{2} + \delta \right|^l u^p = 0, \quad 0<t<1, \\[3ex]
   u(0)=u(1)=0.
  \end{array} 
  \right.
  \label{fracHenon2}
 \end{align}
 has at least three positive nondegenerate solutions for almost everywhere
 $\delta\in(-\delta_0,\delta_0)$.
\end{lemma}

\begin{proof}[Proof of Theorem \ref{nonuniquenessexample}]
 By Lemma \ref{nonuniqueness2}, there exists $\delta_0>0$ such that 
 problem \eqref{fracHenon} with $\alpha=2$ has at least three positive nondegenerate 
 solutions for a.e. $\delta\in(-\delta_0,\delta_0)$.
 Hence, Theorem \ref{existencealphanear2} implies that, for each a.e. 
 $\delta\in(-\delta_0,\delta_0)$, there exists $\alpha_1\in[1,2)$ such that
 problem \eqref{fracHenon} has at least three positive solutions
 for $\alpha\in(\alpha_1,2]$. This completes the proof.
\end{proof} 

\section{H\'{e}non equation} \label{sec:11}

\setcounter{section}{11} \setcounter{equation}{0} 

In this section, we provide a proof of Lemma \ref{nonuniqueness2}. To this end, we consider the one-dimensional  H\'{e}non equation:
\begin{align}\label{Henon}
 \left\{
  \begin{array}{l}
   u'' + |x|^l u^p = 0, \quad -1<x<\zeta, \\[1ex]
   u(-1)=u(\zeta)=0,
  \end{array} 
 \right.
\end{align}
where $l>1$, $p>1$ and $\zeta>-1$.

The number of negative eigenvalues $\mu$ of the following problem \eqref{Morse}
is said to be the {\it Morse index} of $u$:
\begin{align}\label{Morse}
 \left\{
  \begin{array}{l}
   -\phi'' - |x|^l p |u|^{p-1} \phi = \mu \phi, \quad -1<x<\zeta, \\[1ex]
    \phi(-1)=\phi(\zeta)=0.
  \end{array}
 \right.
\end{align}
We have the following result has been obtained in   
\cite[Theorem 1.3]{T2} and \cite[Proposition 1.1 and Lemma 4.1]{ST}.

\begin{proposition}\label{U}
 Let $l>1$ and $p>1$ satisfy $(p-1)l\ge 4$ and let $\zeta=1$.
 Then problem \eqref{Henon} has a unique positive even solution $U$, whose 
 Morse index  is $2$ and $U$ is nondegenerate.
\end{proposition}

We also have the following result from \cite[Proposition 2.3]{T2}.

\begin{proposition}\label{Morse=zero}
 Let $l>1$, $p>1$ and $\zeta>-1$.
 Let $u$ be a positive solution of \eqref{Henon}, and let $w$ be a 
 unique solution of the initial value problem
 \begin{align}\label{LinearMorse}
  \left\{
  \begin{array}{l}
   w'' + |x|^l p |u|^{p-1} w = 0, \quad -1<x\le\zeta, \\[1ex]
   w(-1)=0, \quad w'(-1)=1.
  \end{array}
  \right.
 \end{align}
 Then the Morse index of $u$ is equal to the number of zeros of $w$ in 
 $(-1,\zeta)$.
 Moreover, $u$ is nondegenerate if and only if $w(\zeta)\ne0$.
\end{proposition}

For each $\beta>0$, let $u(x,\beta)$ denote a unique solution of 
the initial value problem
\begin{align}
 \left\{
  \begin{array}{l}
   u'' + |x|^l |u|^{p-1} u = 0, \quad x>-1, \\[1ex]
   u(-1)=0, \quad u'(-1)=\beta.
  \end{array}
 \right.
 \label{initial}
\end{align}
From a general theory on ordinary differential equations and 
a standard argument, we conclude that the solution $u(x,\beta)$ exists on 
$[-1,\infty)$, is unique, $u(x,\beta)$, $u'(x,\beta)$ are $C^1$ 
functions on the set $[-1,\infty)\times(0,\infty)$, and 
$\frac{\partial}{\partial\beta}u(x,\beta)$ is a solution of \eqref{LinearMorse} with
$u=u(x,\beta)$.
Lemma 2.1 in \cite{T2} implies that $u(x,\beta)$ has a zero in $(-1,\infty)$
for each $\beta>0$.
Let $z(\beta)$ denote the smallest zero of $u(x,\beta)$ in $(-1,\infty)$.
Since $u(x,\beta)>0$ for $(-1,z(\beta))$, from the uniqueness of the initial 
value problem, it follows that $u'(z(\beta),\beta)<0$.
Since $u(z(\beta),\beta)=0$, by the implicit function theorem, 
we deduce that $z\in C^1(0,\infty)$ and
\begin{align}
 z'(\beta) 
 = -\frac{\frac{\partial}{\partial\beta}u(z(\beta),\beta)}
         {u'(z(\beta),\beta)}, \quad \beta>0.
 \label{z'}
\end{align}

\begin{lemma}\label{z'>0}
 Let $l>1$ and $p>1$ satisfy $(p-1)l\ge 4$.
 Then there exists $\beta_0>0$ such that $z(\beta_0)=1$ and $z'(\beta_0)>0$.
\end{lemma}

\begin{proof} 
 We use Proposition \ref{U}.
 Then problem \eqref{Henon} with $\zeta=1$ has a unique positive even 
 solution $U$,
 whose Morse index  is $2$, and $U$ is nondegenerate.
 We set $\beta_0=U'(-1)$.
 Then $z(\beta_0)=1$ and $U(x)=u(x,\beta_0)$ for $x\in[-1,1]$.
 Let $w$ be a unique solution of \eqref{LinearMorse} with $u=U$.
 From Proposition \ref{Morse=zero}, it follows that $w$ has exact two zeros in
 $(-1,1)$ and $w(1)\ne0$.
 Thus, $w(1)>0$.
 Since $\frac{\partial}{\partial\beta}u(x,\beta_0)$ is a solution of 
 \eqref{LinearMorse} with $u=u(x,\beta_0)$, we conclude that 
 $w(x)=\frac{\partial}{\partial\beta}u(x,\beta_0)$ for $x\in[-1,1]$.
 By \eqref{z'}, we get $z'(\beta_0)>0$.
\end{proof} 

Applying Lemmas 3.1 and 3.3 in \cite{T2}, we get the following result.

\begin{lemma}\label{z}
 Let $l>1$ and $p>1$.
 Then there exist $\beta_*$ and $\beta^*$ such that $0<\beta_*<\beta^*<\infty$,
 $z(\beta)>1$ for $0<\beta\le\beta_*$ and $z(\beta)<1$ for $\beta\ge\beta^*$.
\end{lemma}

\begin{proof}[Proof of Lemma \ref{nonuniqueness2}]
 By Lemmas \ref{z'>0} and \ref{z}, the intermediate value theorem shows that
 there exists $\varepsilon_0\in(0,1)$ such that, 
 for each $\zeta\in(1-\varepsilon_0,1+\varepsilon_0)$,
 there exist $\beta_1(\zeta)$, $\beta_2(\zeta)$, $\beta_3(\zeta)$ satisfying
 $0<\beta_1(\zeta)<\beta_2(\zeta)<\beta_3(\zeta)<\infty$ and
 $z(\beta_1(\zeta))=z(\beta_2(\zeta))=z(\beta_3(\zeta))=\zeta$.
 This means that, for each $\zeta\in(1-\varepsilon_0,1+\varepsilon_0)$, 
 $u(x,\beta_i(\zeta))$, $i=1$, $2$, $3$ are positive solutions of \eqref{Henon}.
 By Sard's lemma, we find that 
 $z'(\beta_i(\zeta))\ne 0$ for $i=1$, $2$, $3$ and a.e.
 $\zeta\in(1-\varepsilon_0,1+\varepsilon_0)$.
 Hence, by \eqref{z'}, we have
 \begin{align*}
  \frac{\partial}{\partial\beta}u(\zeta,\beta_i(\zeta)) \ne 0,
  \quad i=1,2,3
 \end{align*}
 for a.e. $\zeta\in(1-\varepsilon_0,1+\varepsilon_0)$.
 Since $\frac{\partial}{\partial\beta}u(x,\beta_i(\zeta))$, $i=1$, $2$, $3$ 
 are solutions of \eqref{LinearMorse} with $u=u(x,\beta_i(\zeta))$, 
 Proposition \ref{Morse=zero} implies that $u(x,\beta_i(\zeta))$, 
 $i=1$, $2$, $3$ are nondegenerate 
 for a.e. $\zeta\in(1-\varepsilon_0,1+\varepsilon_0)$.
 We set
 \begin{align*}
  v_i(t)= (1+\zeta)^\frac{l+2}{p} u((1+\zeta)t-1,\beta_i(\zeta)).
 \end{align*}
 Then we conclude that, for a.e. $\zeta\in(1-\varepsilon_0,1+\varepsilon_0)$,
 $v_1$, $v_2$ and $v_3$ are positive nondegenerate solutions of problem 
 \begin{align*}
  \left\{
  \begin{array}{l}
   v'' + \left| t -\frac{1}{2} + \frac{\zeta-1}{2(1+\zeta)} \right|^l v^p = 0,
   \quad 0<t<1, \\[2ex]
   v(0)=v(1)=0.
  \end{array} 
  \right.
 \end{align*}
 This completes the proof.
\end{proof} 

\noindent{\bf Acknowledgements} \
The authors would like to thank the reviewers for helpful comments on
important revisions that improve the results of our manuscript.


\end{document}